\documentclass[10pt,twocolumn,twoside,journal]{IEEEtran} 

\usepackage{cite}
\usepackage{amsmath} 
\usepackage{amssymb}  
\usepackage{amsfonts}
\usepackage{mathrsfs}
\usepackage{mathtools}
\usepackage[dvipsnames]{xcolor}
\usepackage[thmmarks, amsmath]{ntheorem}
\usepackage{enumitem,xcolor}
\usepackage{circuitikz}
\usepackage{pgfplots}
\usepackage{booktabs}
\usepackage[bbgreekl]{mathbbol}
\usepackage{colortbl}
\pgfplotsset{compat=1.14}

\DeclareSymbolFontAlphabet{\mathbb}{AMSb}
\DeclareSymbolFontAlphabet{\mathbbl}{bbold}

\DeclareMathOperator{\diag}{diag}

\DeclareMathOperator{\diff}{d}

\DeclareMathOperator{\rank}{rank}

\providecommand{\norm}[1]{\lVert#1\rVert}

\newcommand{\R}{{\mathbb R}}
\newcommand{\mc}{\mathcal}
\newcommand{\ddt}{\tfrac{\diff}{\diff \! t}}

\newtheorem{theorem}{Theorem}
\newtheorem{lemma}{Lemma}
\newtheorem{proposition}{Proposition}
\newtheorem{assumption}{Assumption}

\newtheorem{definition}{Definition}
\newtheorem{condition}{Condition}

\graphicspath{{Figures/} }
\allowdisplaybreaks

\begin{document}

\title{\LARGE \bf  The effect of transmission-line dynamics on  grid-forming dispatchable virtual oscillator control  \thanks{This work was partially funded by the European Union's Horizon 2020 research and innovation programme under grant agreement N$^\circ$ 691800. This article reflects only the authors' views and the European Commission is not responsible for any use that may be made of the information it contains.}}
\author{Dominic Gro\ss{}, Marcello Colombino, Jean-S\'ebastien Brouillon,  and  Florian D\"orfler  \thanks{D. Gro\ss{}, J-S Brouillon and F. D\"orfler, are with the Automatic Control Laboratory, ETH Z\"urich, Switzerland, M. Colombino is with the Electrical and Computer Engineering Department at McGill University, Canada; e-mail:\{grodo,jeanb,dorfler\}@ethz.ch,  marcello.colombino@mcgill.ca}}

\maketitle

\begin{abstract}
In this work, we analyze the effect of transmission line dynamics on grid-forming control for inverter-based AC power systems. In particular, we investigate a dispatchable virtual oscillator control (dVOC) strategy that was recently proposed in the literature. When the dynamics of the transmission lines are neglected, i.e., if an algebraic model of the transmission network is used, dVOC ensures almost global asymptotic stability of a network of AC power inverters with respect to a pre-specified solution of the AC power-flow equations. While this approximation is typically justified for conventional power systems, the electromagnetic transients of the transmission lines can compromise the stability of an inverter-based power system. In this work, we establish explicit bounds on the controller set-points, branch powers, and control gains that guarantee almost global asymptotic stability of dVOC in combination with a dynamic model of the transmission network.
\end{abstract}

\section{Introduction}
The electric power grid is undergoing a period of unprecedented change. A major transition is the replacement of bulk generation based on synchronous machines by renewable generation interfaced via power electronics. This gives rise to scenarios in which either parts of the transmission grid or an islanded distribution grid may operate without conventional synchronous generation. In either case,  the loss of synchronous machines poses a great challenge because today's power system operation heavily relies on their self-synchronizing dynamics, rotational inertia, and resilient controls.

The problem of synchronization of grid-forming power inverters has been widely studied in the recent literature. A grid-forming inverter is not limited to power tracking but acts as a controlled voltage source that can change its power output (thanks to storage or curtailment), and is controlled to contribute to the stability of the grid. Most of the common approaches of grid-forming control focus on droop control~\cite{MCC-DMD-RA:93,Guerrero2015,simpson2017voltage}. Other popular approaches are based on mimicking the physical characteristics and controls of synchronous machines~\cite{zhong2011synchronverters,JAD18,SDA-SJA:13} or controlling inverters to behave like virtual Li\'enard-type oscillators~\cite{johnson2014synchronization,MS-FD-BJ-SD:14b,LABT-JPH-JM:13,johnson2016synthesizing}. While strategies based on machine-emulation are compatible with the legacy power system, they use a system (the inverter) with fast actuation but almost no inherent energy storage  to mimic a system (the generator) with slow actuation but significant energy storage (in the rotating mass) and may also be ineffective due to deteriorating effects of time-delays in the control loops and current limits of inverters \cite{ENTSOE16}. While some form of power curtailment or energy storage is essential to maintain grid stability, it is not clear that emulating a machine behavior is the preferable option. Virtual oscillator control (VOC) is a promising approach because it can globally synchronize an inverter-based power system. However, {the nominal power injection cannot be specified in the original VOC approach~\cite{johnson2014synchronization,MS-FD-BJ-SD:14b,LABT-JPH-JM:13,johnson2016synthesizing}, i.e., it cannot be dispatched.} Likewise, all theoretic investigations are  limited to synchronization with identical angles and voltage magnitudes. {For passive loads it can be shown that power is delivered to the loads \cite{LABT-JPH-JM:13}, but the power sharing by the inverters and their voltage magnitudes are determined by the load and network parameters.}
Finally, the authors recently proposed a dispatchable virtual oscillator control (dVOC) strategy~\cite{colombino2017global,colombino2017global2}, which relies on synchronizing harmonic oscillators through the transmission network, ensures almost global asymptotic stability of an inverter-based AC power systems with respect to a desired solution of the AC power-flow equations{, and has been experimentally validated in \cite{GC+19}}. 

In order to simplify the analysis, the dynamic nature of transmission lines is typically neglected in the study of power system transient stability and synchronization. In most of the aforementioned studies, an algebraic model of the transmission network is used, i.e., the relationship between currents and voltages is modeled by the admittance matrix of the network. 
This approximation is justified in a traditional power network, where bulk generation is provided by synchronous machines with very slow time constants (seconds) and can be made rigorous using time-scale separation arguments~\cite{curi2017control}. As power inverters can be controlled at much faster time-scales~(milliseconds), when more and more synchronous machines are replaced by inverter-based generation, the transmission line dynamics, which are typically not accounted for in the theoretical analysis, can compromise the stability of the power-network. For droop controlled microgrids this phenomenon has been noted in~\cite{Guerrero2015,vorobev2017high} and it can be verified experimentally for all control methods listed above. Moreover, in \cite{vorobev2017high,vorobev2017cdc} explicit and insightful bounds on the control gains are obtained via small signal stability analysis for a steady state with zero relative angle.

In this work, we study the dVOC proposed in~\cite{colombino2017global}, which renders power inverters interconnected through an algebraic network model almost globally asymptotically stable \cite{colombino2017global2}. We provide a Lyapunov characterization of almost global stability. This technical contribution is combined with ideas inspired by singular perturbation theory \cite{PKC82,K02} to construct a Lyapunov function and an explicit analytic stability condition that guarantees almost global asymptotic stability of the full power-system including transmission line dynamics. {Broadly speaking, this stability condition makes the time-scale separation argument rigorous by quantifying how large the time-scale separation between the inverters and the network needs to be to ensure stability.
In addition to capturing the influence of the network dynamics our stability condition is not restricted to operating points with zero relative angles, accounts for the network topology, and provides explicit bounds the control gains and set-points for the power injections.} In particular, the stability condition links the achievable power transfer and maximum control gain to the connectivity of the graph of the transmission network. Moreover, we show that the only undesirable steady-state of the closed loop, which corresponds to voltage collapse, is exponentially unstable. By guaranteeing almost global stability of dVOC our results provide a theoretical foundation  for using dVOC for non-standard operating conditions such as black starts and islanded operation, and alleviate the challenging problem of characterizing the region of attraction of networks of droop-controlled inverters. We use a simplified three-bus transmission system to show that our stability condition is not overly conservative and, to illustrate the results, we provide simulations of the IEEE 9-bus system.

In Section~\ref{sec.model} we introduce the problem setup. Section~\ref{sec.lyap} provides definitions and a preliminary technical results. The main result is given in Section~\ref{sec.sing.pert}, Section~\ref{sec.example} presents numerical examples, and Section~\ref{sec.conclusion} provides the conclusions.

\subsection*{Notation}
The set of real numbers is denoted by $\R$, $\R_{\geq 0}$ denotes the set of nonnegative reals, and $\mathbb N$ denotes the natural numbers. Given ${\theta\in[-\pi,\pi]}$ the 2D rotation matrix is given by
\begin{align*}
R(\theta) \coloneqq \begin{bmatrix}\cos(\theta) & -\sin(\theta)\\  \sin(\theta) & \cos(\theta) \end{bmatrix} \in \R^{2 \times 2}.
\end{align*}
Moreover, we define the $90^\circ$ rotation matrix $J \coloneqq R(\pi/2)$ that can be interpreted as an embedding of the complex imaginary unit $\sqrt{-1}$ into $\R^2$. Given a matrix $A$, $A^\mathsf{T}$ denotes its transpose. We use $\|A\|$ to indicate the induced 2-norm of $A$. We write $A\succcurlyeq0$  $(A\succ0)$ to denote that $A$ is symmetric and positive semidefinite (definite). For column vectors $x\in\R^n$ and $y\in\R^m$ we use $(x,y) = [x^\mathsf{T}, y^\mathsf{T}]^\mathsf{T} \in \R^{n+m}$ to denote a stacked vector and $\norm{x}$ denotes the Euclidean norm. {The absolute value of a scalar $y \in \R$ is denoted by $|y|$.} Furthermore, $I_n$ denotes the identity matrix of dimension $n$, and $\otimes$ denotes the Kronecker product. Matrices of zeros of dimension $n \times m$ are denoted by $\mathbbl{0}_{n\times m}$ and $\mathbbl{0}_{n}$ denotes column vector of zeros of length $n$. We use $\norm{x}_C \coloneqq \min_{z \in \mathcal{C}} \norm{z-x}$ to denote the distance of a point $x$ to a set $\mathcal{C}$. {Moreover, the cardinality of a set $\mc E$ is denoted by $|\mc E|$.} We use $\varphi_f(t,x_0)$ to denote the solution of $\ddt x = f(x)$ at time $t \geq 0$ starting from the initial condition $x(0)=x_0$ at time $t_0=0$.

\section{Modeling and control of an inverter-based power network}\label{sec.model}
In this section, we introduce the model of the inverter-based power grid that will be studied throughout the paper.

\subsection{Inverter-based power grid}
We study the control of $N$ three-phase inverters interconnected by a resistive-inductive network. All electrical quantities in the network are assumed to be balanced. This allows us to work in $\alpha\beta$ coordinates obtained by applying the well-known Clarke transformation to the three-phase variables~\cite{clarke1943circuit}. We model each inverter as a controllable voltage source. We model the transmission network as a simple graph (i.e., an undirected graph containing no self-loops or multiple edges) denoted by $\mc G~=~(\mc N, \mc E,\mc W)$, where $\mc N = \{1,..., N\}$ is the set of nodes corresponding to the inverters, $\mc W$ is the set of edge weights, and $\mc E \subseteq (\mc N \times \mc N) \setminus \cup_{i \in \mc N} (i,i)$, with $|\mc E|=M$, is the set of edges corresponding to the transmission lines. To each inverter we associate a terminal voltage $\underline{v}_k\in\R^2$, that can be fully controlled, and an output current $\underline{i}_{o,k}\in\R^2$ flowing out of the inverter and into the network. To each transmission line, we associate a current $\underline{i}_l \in \R^2$, an impedance matrix $Z_l \coloneqq I_2 r_{l} + \omega_0 J \ell_{l}$, and an admittance matrix $Y_l \coloneqq Z^{-1}_l$ with
\begin{align}\label{eq.weights}
\norm{Y_l} = \frac{1}{\sqrt{r_l^2+\omega_0^2\ell_l^2}},
\end{align}
where $r_l \in \mathbb{R}_{>0}$ and $\ell_l \in \mathbb{R}_{>0}$ are the resistance and inductance of the line $l \in \{1,\ldots,M\}$ respectively, and $\omega_0 \in \mathbb{R}_{\geq 0}$ is the nominal operating frequency of the power system. Moreover, we define the edge weights $\mc W = \{\norm{Y_l}\}_{l=1}^M$ of the graph. The oriented incidence matrix of the graph is denoted by $B \in \{-1,0,1\}^{N \times M}$ and, duplicating each edge for the $\alpha$ and $\beta$ components, we define $\mc B \coloneqq B\otimes I_2$. The dynamics of the transmission line currents $\underline{i} \coloneqq (\underline{i}_1,\ldots,\underline{i}_M) \in \R^{2M}$ are given by
\begin{align}\label{eq.network.equations}
L_T \ddt \underline{i} & = -R_T \underline{i} + \mc B^\mathsf{T} \underline{v},
\end{align}
where $R_T\coloneqq \diag(\{r_l\}_{l=1}^M)  \otimes I_2$ is the matrix of line resistances $r_l \in \mathbb{R}_{>0}$, $L_T \coloneqq \diag(\{\ell_l\}_{l=1}^M) \otimes I_2$ is the matrix of line inductances $\ell_l \in \mathbb{R}_{>0}$, and $\underline{v} \coloneqq (\underline{v}_1,\ldots,\underline{v}_N) \in \R^{2N}$ collects the terminal voltages of all inverters. We define the vector of inverter output currents as $\underline{i}_o \coloneqq (\underline{i}_{o,1},\ldots,\underline{i}_{o,N}) \in \R^{2N}$ given by $\underline{i}_o =   \mc B \, \underline{i}$. For convenience we use $\ell_{jk}$, $r_{jk}$, and $Y_{jk}$, to denote the inductance, resistance, and impedance of the line connecting node $j$ and node $k$.

\subsection{Quasi-steady-state network model}
To obtain the \emph{quasi-steady-state} approximation (also known as \emph{phasor approximation}) of~\eqref{eq.network.equations} we perform the following change of variables to a rotating reference frame:
\begin{align}\label{eq.rotframe}
\underline{v} = \mc R(\omega_0 t) v, \quad \underline{i}_o =  \mc R(\omega_0 t) i_o, \quad \underline{i} = I_M \otimes R(\omega_0 t) i, 
\end{align}
where $\mc R(\theta) \coloneqq I_N \otimes R(\theta)$. Using $\mc J_M = I_M \otimes J$ and the impedance matrix $Z_T \coloneqq R_T + \omega_0 \mc J_M L_T$, the dynamics~\eqref{eq.network.equations} in the rotating frame become 
\begin{align}\label{eq.network.rot}
L_T \ddt i  = -Z_T i +  \mc B^\mathsf{T} v.
\end{align}
The quasi-steady-state model is obtained by considering the steady-state map of the exponentially stable dynamics \eqref{eq.network.rot}.
\begin{definition}{\bf(Quasi-steady-state network model)}\label{def:qs}\\
The quasi-steady-state model of the transmission line currents $i$ and output currents $i_o$ is given by the steady-state map $i^s(v) \coloneqq Z^{-1}_T \mc B^\mathsf{T} v$ of \eqref{eq.network.rot} and $i^s_o(v) \coloneqq 
\mc Y v$ with $\mc Y \coloneqq \mc B Z^{-1}_T \mc B^\mathsf{T}$.
\end{definition}
In conventional power systems the approximation $i = i^s(v)$ is typically justified due to the pronounced time-scale separation between the dynamics of the transmission lines and the synchronous machine dynamics. Therefore, the dynamic nature of the transmission lines is typically neglected a priori in transient stability analysis of power systems. However, for inverter-based power systems the electromagnetic transients of the lines have a significant influence on the stability boundaries, and the approximation is no-longer valid \cite{vorobev2017high,vorobev2017cdc}. In this work, we make the time-scale separation argument rigorous and explicit by quantifying how large the time-scale separation needs to be to ensure stability. In order to prove the main result of this manuscript, we require the following standing assumption.
\begin{assumption}{\bf(Uniform inductance-resistance ratio)}\label{ass.constant.ratio}\\
The ratio between the inductance and resistance of every transmission line is constant, i.e., for all $(j,k) \in\mc E$ it holds that $\frac{\ell_{jk}}{r_{jk}} = \rho \in \R_{>0}$.
\end{assumption}
Assumption~\ref{ass.constant.ratio} is typically  approximately satisfied for transmission lines at the same voltage level, but not satisfied across different voltage levels. We note that the main salient features uncovered by our theoretical analysis can still be observed even when Assumption~\ref{ass.constant.ratio} only holds approximately (see Section \ref{sec:example:ieee9}). 
Finally, let us define the angle $\kappa  \coloneqq \tan^{-1}( \rho \omega_0)$, the graph Laplacian $L \coloneqq B\, \diag(\{\norm{Y_l}\}_{l=1}^M) B^\mathsf{T}$, and the \emph{extended} Laplacian $\mc L \coloneqq L \otimes I_2$. Note that under Assumption \ref{ass.constant.ratio}, it can be verified that $\mc R(\kappa) i^s_o(v) = \mc R(\kappa) \mc Y v = \mc L v$.

\subsection{Control objectives}\label{sec.control.obj}
In this section we formally specify the \emph{system-level} control objectives for an inverter-based AC power system that need to be achieved with \emph{local} controllers. We begin by defining instantaneous active and reactive power.
\begin{definition}{\bf(Instantaneous Power)}\label{def:IP}\\
Given the voltage $\underline{v}_k$ at node $k \in \mc N$ and the output current $\underline{i}_{o,k}$, we define the \emph{instantaneous active power} $p_k \coloneqq \underline{v}_{k}^\mathsf{T}  \underline{i}_{o,k} \in \mathbb{R}$ and the \emph{instantaneous reactive power} $q_k \coloneqq \underline{v}_k^\mathsf{T} J \underline{i}_{o,k} \in \mathbb{R}$. Moreover, for all $(j,k) \in \mc E$, we define the \emph{instantaneous} active and reactive \emph{branch powers} $p_{jk} \!\coloneqq \underline{v}^\mathsf{T}_k \underline{i}_{jk}$ and $q_{jk} \!\coloneqq \underline{v}^\mathsf{T}_k J \underline{i}_{jk}$.
\end{definition}
Active and reactive power injections cannot be prescribed arbitrarily to each inverter in a network but they need to be consistent with the power-flow equations~\cite{kundur1994power}. To this end, we introduce steady-state voltage angles $\theta^\star_{k1}$ relative to the node with index $k=1$, and we define the relative steady-state angles $\theta^\star_{jk} \coloneqq \theta^\star_{j1} - \theta^\star_{k1}$ for all $(i,j) \in \mc N \times \mc N$.
\begin{condition}\textbf{(Consistent set-points)}\label{cond.consistent}\\ 
The set-points $p_k^\star \in \mathbb{R}$, $q_k^\star \in \mathbb{R}$, $v_k^\star \in \mathbb{R}_{>0}$ for active power, reactive power, and voltage magnitude respectively, are consistent with the power flow equations, i.e., for all $(j,k) \in \mc E$ there exist relative angles ${\theta^\star_{jk}} \in (-\pi,\pi]$ and steady-state \emph{branch powers} $p^\star_{jk} \in \mathbb{R}$ and $q^\star_{jk} \in \mathbb{R}$ given by
\begin{align*}
 p^\star_{jk} \!\!&\coloneqq\! \norm{Y_{jk}}^2 \!\big(v^{\star2}_k r_{jk} - v^{\star}_k v^{\star}_j
( r_{jk}\cos(\theta_{jk}^\star )\!+\omega_0\ell_{jk} \sin(\theta_{jk}^\star ))\big),\\
 q^\star_{jk} \!\!&\coloneqq\! \norm{Y_{jk}}^2 \!\big( v_k^{\star2} \omega_0\ell_{jk} \!-\! v_k^{\star}v_j^{\star}(  \omega_0\ell_{jk}\!\cos(\theta_{jk}^\star )\!-r_{jk} \sin(\theta_{jk}^\star ))\big),
\end{align*}
such that $p^\star_k = \sum_{(j,k)\in\mc E} p^\star_{jk}$ and $q^\star_k = \sum_{(j,k)\in\mc E} q^\star_{jk}$  holds for all $k \in \mc N$.
\end{condition}
Given a set of \emph{consistent} power, voltage and frequency set-pints (according to Condition~\ref{cond.consistent}), the objective of this paper is to design a decentralized controller of the form
\begin{align}\label{eq.voltage.dyn.ol}
\ddt \underline{v}_k & =u_k(\underline{v}_k,\underline{i}_{o,k}),
\end{align}
that induces the following steady-state behavior:
\begin{itemize}
\item \textbf{Synchronous frequency:} Given a desired synchronous frequency $\omega_0 \in \R_{\ge0}$, at steady state it holds that:
\begin{align}\label{eq.obj.sync}
\ddt {\underline{v}_k} = \omega_0J\,\underline{v}_k,\!\quad {k}~\!\in~\mc N;
\end{align}
\item \textbf{Voltage magnitude:} Given voltage magnitude set-points $v^\star_k \in \R_{>0}$, at steady state it holds that:
\begin{align}\label{eq.obj.norm}
\norm{\underline{v}_k}= v^\star_k, \quad {k}~\!\in~\mc N;
\end{align}
\item \textbf{Steady-state currents:} it holds that $\ddt \underline{i} = \omega_0 \mathcal{J}_M \underline{i}$;
\item \textbf{Power injection:} At steady state, each inverter injects the prescribed active and reactive power, i.e.,
\begin{align}\label{eq.obj.power}
\underline{v}_k^\mathsf{T} \underline{i}_{o,k} =p_k^\star, \quad \underline{v}_k^\mathsf{T} J \underline{i}_{o,k}=q_k^\star,
\quad {k}~\!\in~\mc N.
\end{align}
\end{itemize}

Equation~\eqref{eq.obj.sync} specifies that, at steady state, all voltages in the power network evolve as sinusoidal signals with frequency $\omega_0$. Equations~\eqref{eq.obj.norm} and~\eqref{eq.obj.power} specify the voltage magnitude and the power injection at every node. The \emph{local nonlinear} specification \eqref{eq.obj.power} on the steady-state power injection can be equivalently expressed as \emph{non-local linear} specification on the voltages as follows (cf. \eqref{eq.powers.ref} and \eqref{eq.powang}):
\begin{itemize}
\item \textbf{Phase locking:} Given relative angles {$\theta^\star_{k1} \in [-\frac{\pi}{2},\frac{\pi}{2}]$}, at steady state it holds that:
\begin{align}\label{eq.obj.theta}
\frac{\underline{v}_k}{v^\star_k} -  R({\theta_{k1}^\star})\, \frac{\underline{v}_1}{v^\star_1} = 0, \quad k~\!\in~\mc N\backslash\{1\}.
\end{align}
\end{itemize}
\subsection{Dispatchable virtual oscillator control}\label{sec.control.cl.design}
For every inverter $k \in \mc N$, we define the dispatchable virtual oscillator control (dVOC) 
\begin{align}\label{eq.control.law}
 u_k \coloneqq \omega_0 J \underline{v}_k + \eta\left[K_k \underline{v}_k - R(\kappa)\underline{i}_{o,k} + \alpha \Phi_k(\underline{v}_k)\underline{v}_k\right],
\end{align}
with gains $\eta \in \mathbb{R}_{>0}$, $\alpha \in \mathbb{R}_{>0}$, and 
\begin{align*}
K_k = \frac{1}{v_k^{\star2}} R(\kappa) \begin{bmatrix} p_k^\star & q_k^\star\\ -q_k^\star & p_k^\star
\end{bmatrix}, \quad
\Phi_k(\underline{v}_k)  \coloneqq   \frac{v_k^{\star2}-\norm{\underline{v}_k}^2}{v_k^{\star2}} I_2.
\end{align*}
Note that the term $\omega_0 J \underline{v}_k$ induces a harmonic oscillation of $\underline{v}_k$ at the nominal frequency $\omega_0$. Moreover, $\Phi_k(\underline{v}_k)\underline{v}_k$ can be interpreted as a voltage regulator, i.e.,
depending on the sign of the normalized quadratic voltage error $\Phi_k(\underline{v}_k)$ the voltage vector $\underline{v}_k$ is scaled up or down. Finally, the term $K_k \underline{v}_k - R(\kappa)\underline{i}_{o,k}$ can be interpreted either in terms of {tracking power set-points (i.e., \eqref{eq.obj.power}), or in terms of phase synchronization (i.e., \eqref{eq.obj.theta}).} To establish the connection to {tracking power set-points}, let $e_{p,k}$ and $e_{q,k}$ denote the weighted difference between the actual and desired power injections defined by ${v^{\star2}_k} e_{p,k} \coloneqq \norm{\underline{v}_k}^2 p^\star_k -  {v^{\star2}_k} p_k$ and ${v^{\star2}_k} e_{q,k} \coloneqq \norm{\underline{v}_k}^2 q^\star_k - {v^{\star2}_k} q_k$. A straightforward algebraic manipulation reveals that
\begin{align}\label{eq.pspecs.rot}
\frac{1}{\norm{\underline{v}_k}^2} \begin{bmatrix} \underline{v}_k  & J \underline{v}_k \end{bmatrix}
 R(\kappa) \begin{bmatrix} e_{p,k} \\ -e_{q,k} \end{bmatrix}\! = 
K_k \underline{v}_k - R(\kappa)\underline{i}_{o,k},
\end{align}
i.e., normalizing the vector $(e_{p,k},-e_{q,k})$ and transforming it into a rotating frame attached to $\underline{v}_k$ results in the synchronizing dVOC term. In particular, in a purely inductive grid it holds that $R(\kappa) = J$ and $e_{p,k}$ corresponds to a component that is tangential to $\underline{v}_k$, i.e., its rotational speed, and $e_{q,k}$ corresponds to a component that is radial to $\underline{v}_k$, i.e., its change in magnitude. In other words, in this case, dVOC resembles standard droop control. For networks which are not purely inductive, the rotation $R(\kappa)$ results in a mixed nonlinear droop behavior. 
By construction of the control law~\eqref{eq.control.law}, if $p_k=p^\star_k$, $q_k = q_k^\star$, and $\|\underline v_k\|=v_k^\star$, then $\underline{u}_k = \omega_0 J \underline{v}_k$, leaving all voltage vectors to rotate synchronously at frequency $\omega_0$.
{To establish the connection to phase synchronization (i.e., \eqref{eq.obj.theta}), let $e_\theta(\underline{v}) \coloneqq (e_{\theta,1}(\underline{v}),\ldots,e_{\theta,N}(\underline{v}))$ denote the vector of weighted phase errors defined by
\begin{align}\label{eq.distrib.feedback}
 e_{\theta,k}(\underline{v}) \coloneqq \sum\nolimits_{(j,k)\in\mc E} \norm{Y_{jk}} (\underline{v}_j - \tfrac{v^\star_j}{v^\star_k}R(\theta_{jk}^\star) \underline{v}_k).
\end{align}
Under the \emph{quasi-steady-state} approximation $\underline{i}_{o,k} = \underline{i}^s_o(\underline{v}) = i^s_o(\underline{v})$ the \emph{local} feedback in~\eqref{eq.control.law} is identical to a \emph{distributed} synchronizing feedback of the weighted phase errors~\eqref{eq.distrib.feedback}. This is made precise in the following proposition obtained by combining Proposition 1 and Proposition 2 of \cite{colombino2017global}.
\begin{proposition}{\bf(Synchronizing feedback)}\label{prop.sync}\\
 Consider set-points $p^\star_k$, $q^\star_k$, $v^\star_k$, and steady-state angles $\theta^\star_{jk}$  that satisfy Condition \ref{cond.consistent}. It holds that $K_k\underline{v}_k - R(\kappa)\underline{i}^s_{o,k} = e_{\theta,k}(\underline{v})$.
\end{proposition}
The proof is given in the Appendix.
Proposition \ref{prop.sync} highlights how dVOC infers information on the phase difference of the voltages $\underline{v}_k$ from the currents $i_o(\underline{v})$.} In~\cite{colombino2017global2} the authors use this fact to prove stability of the controller~\eqref{eq.control.law} under the quasi-steady-state approximation $\underline{i}_o = i^s_o(\underline{v})$. However, the line dynamics~\eqref{eq.network.equations} ``delay" the propagation of this information and give rise to stability concerns that motivate this work.

\section{Almost global asymptotic stability of set-valued control specifications}\label{sec.lyap}
In order to state the main result of the paper, we require the following technical definitions and preliminary results, which are used in Section~\ref{sec.sing.pert} to analyze the stability properties of the inverter-based power system. 
We begin by defining almost global asymptotic stability (see \cite{A04}) with respect to a set $\mc C$.
\begin{definition}{\bf(Almost global asymptotic stability)}\label{def:ags}
A dynamic system $\ddt x = f(x)$ is called almost globally asymptotically stable with respect to a compact set $\mc C$ if 
\begin{enumerate}[label=(\roman*)]
 \item it is almost globally attractive with respect to $\mc C$, i.e.,
  \begin{align}
     \lim_{t\to\infty} \norm{\varphi_f(t,x_0)}_{\mc C} = 0
   \end{align}
        holds for all $x_0 \notin \mc Z$ and $\mc Z$ has zero Lebesgue measure, 
        \item it is Lyapunov stable with respect to $\mc C$, i.e., for every $\varepsilon \in \mathbb{R}_{>0}$ there exists $\delta \in \mathbb{R}_{>0}$ such that
  \begin{align}\label{eq:lyapstab}
   \norm{x_0}_{\mc C} < \delta \implies \norm{\varphi_f(t,x_0)}_{\mc C} < \varepsilon, \qquad \forall t \geq 0.
   \end{align}
\end{enumerate}
\end{definition}
Moreover, we define two classes of comparison functions which are used to establish stability properties of the system.
\begin{definition}{\bf(Comparison functions)}
 A function $\chi_c:\mathbb{R}_{\geq0} \to \mathbb{R}_{\geq0}$ is of class $\mathscr{K}$ if it is continuous, 
strictly increasing and $\chi_c(0)=0$; it is
 of class  $\mathscr{K}_{\infty}$ if it is a $\mathscr{K}$-function and $\chi_c(s)\to\infty$ as $s\to\infty$. 
\end{definition}
The following Theorem provides a Lyapunov function characterization of almost global asymptotic stability.
\begin{theorem}{\bf(Lyapunov function)}\label{thm.AGAS}
Consider a Lipschitz continuous function $f: \mathbb{R}^n \to \mathbb{R}^n$, a compact set $\mc C \subset \mathbb{R}^n$, and an invariant set $\mc U \subset \mathbb{R}^n$ {(}i.e., $\varphi_f(t,x_0) \in \mc U$ for all $t \in \R_{>0}$ and all $x_0 \in \mc U${)}, such that $\mc C \cap \mc U = \emptyset$. Moreover, consider a continuously differentiable function $\mc V: \mathbb{R}^n \to \mathbb{R}_{>0}$ and comparison functions $\chi_1, \chi_2 \in \mathscr{K}_\infty$ and $\chi_3 \in \mathscr{K}$ such that 
 \begin{align*}
  \chi_1(\norm{x}_{\mathcal{C}}) \leq \mc V(x) &\leq \chi_2(\norm{x}_{\mathcal{C}})\\
    \ddt  \mc V(x) \coloneqq \frac{\partial \mc V}{\partial x} f(x) &\leq -\chi_3(\norm{x}_{\mc C \cup \, {\mc U}})
  \end{align*} 
  holds for all $x \in \mathbb{R}^n$. Moreover, let 
  \begin{align*}
   \mc Z_{\mc U} \coloneqq \{ x_0 \in \mathbb{R}^n \vert \lim\nolimits_{t\to\infty} \norm{\varphi_{f}(t,x_0)}_{\mc U} = 0 \}.
  \end{align*}
  denote the region of attraction of $\mc U$. If  $\mc Z_{\mc U}$ has zero Lebesgue measure, the dynamics $\ddt x = f(x)$ are almost globally asymptotically stable with respect to $\mathcal{C}$.
\end{theorem}
The proof is given in the Appendix.

\section{Stability analysis of the inverter-based AC power system}\label{sec.sing.pert}
In Section~\ref{sec.model} we introduced the control objective and proposed a control law that admits a fully decentralized implementation. We will now analyze the closed-loop system and provide sufficient conditions for stability. To do so, we will use ideas from singular perturbation analysis to verify the assumptions of Theorem~\ref{thm.AGAS} for the multi-inverter power system. This allows us to prove the main result of the paper: an explicit bound for the control gains and set-points that ensures that the steady-state behavior given by~\eqref{eq.obj.sync}-\eqref{eq.obj.power} is almost globally asymptotically stable according to Definition~\ref{def:ags}. 

\subsection{Dynamics and control objectives in a rotating frame}
In order to simplify the analysis, it is convenient to perform the change of variables \eqref{eq.rotframe} to a rotating reference frame. Defining $\mc K \coloneqq\diag(\{K_k\}_{k=1}^N)$ and $\Phi(\underline{v}) \coloneqq\diag\left(\{\Phi_k(\underline{v}_k)\}_{k=1}^N\right)$ the combination of \eqref{eq.voltage.dyn.ol}, dVOC~\eqref{eq.control.law}, and the line dynamics~\eqref{eq.network.equations} in the new coordinates becomes
\begin{subequations}\label{eq.closed.loop.rot}
\begin{align}
\!\!\ddt v & = \eta(\mc K v - \mc R(\kappa) \mc B\,i + \alpha \Phi(v)v)=:  f_v(v,i), \label{eq.voltage.cl.rot}\\
\!\!\ddt i & = L_T^{-1}\left(-Z_T i +  \mc B^\mathsf{T} v\right) =:  f_i(v,i). \label{eq.lines.cl.rot}\!
\end{align}
\end{subequations}
Moreover, we let  $x \coloneqq (v,i) \in \R^n$ with $n = 2N+2M$, $f(x) \coloneqq (f_v(v,i),f_i(v,i))$, and denote the dynamics of the overall system by $\ddt x = f(x)$. 
To formalize the control objectives \eqref{eq.obj.sync}-\eqref{eq.obj.theta}, we define the sets
\begin{align*}
\mc S &= \left\{v \in \R^{2N} \left\vert\; \frac{v_k}{v^\star_k}=  R(\theta_{k1}^\star)\,\frac{v_1}{v^\star_1},\; \forall k \in \mc N \setminus \{1\} \right.\right\},\\
\mc A &= \left\{v \in \R^{2N} \left\vert\; \norm{v_k} =  v_k^\star,\; \forall k \in \mc N \right.\right\},
\end{align*}
as well as the target set
\begin{align}\label{eq.T}
\mc T\coloneqq \left\{x \in \R^{n} \left\vert\; x = (v,i),~ v \in \mc S\cap\mc A,~ i = i^s(v)  \right.\right\}.
\end{align} 
Note that $i^s(v)$ is the steady-state map of~\eqref{eq.lines.cl.rot}. Note that all elements of $\mc T$ are equilibria for the dynamics in the rotating reference frame~\eqref{eq.closed.loop.rot}. Therefore, in the static frame, they correspond to synchronous sinusoidal trajectories with frequency $\omega_0$. By definition of the sets $\mc S$, $\mc A$, and $\mc T$, they also satisfy all control objectives introduced in Section~\ref{sec.control.obj}. Moreover, we define the union of $\mc T$ and the origin as $\mc T_0 \coloneqq \mc T \cup \{\mathbbl{0}_{n}\}$. 
\subsection{Main result}
We require the following condition to establish almost global asymptotic stability of the power system dynamics~\eqref{eq.closed.loop.rot}.
\begin{condition}{\bf(Stability Condition)}\label{cond.stab}
The set-points $p^\star_k$, $q^\star_k$, $v_k^\star$ and the steady-state angles $\theta^\star_{jk}$ satisfy Condition~\ref{cond.consistent}. There exists a maximal steady-state angle $\bar{\theta}^\star \in [0,\tfrac{\pi}{2}]$ such that $|\theta^\star_{jk}| \! \leq \! \bar{\theta}^\star$ holds for all $(j,k) \! \in \! \mc N \!\times\! \mc N$. For all $k \! \in \! \mc N$, the line admittances $\norm{Y_{jk}}$, the stability margin $c \in \mathbb{R}_{>0}$, the set-points $p^\star_k$, $q^\star_k$, $v_k^\star$, and the gains $\eta \in \mathbb{R}_{>0}$ and $\alpha \in \mathbb{R}_{>0}$ satisfy
\begin{multline*}
\!\!\!\!\!\!\sum_{j : (j,k) \in \mc E} \!\!\!\!\norm{Y_{jk}}\!\left|1\!-\!\frac{v^\star_j}{v^\star_k}\cos(\theta^\star_{jk})\right|+ \alpha {<} \frac{{1\!+\!\cos(\bar{\theta}^\star)}\!}{2} \frac{v^{\star2}_{\min}}{v^{\star2}_{\max}}\lambda_2( L)-c,\\
\eta < \frac{c }{\rho \norm{\mc Y} (c  + 5 \|\mc K -\mc L\|)},
\end{multline*}
where $v^{\star}_{\min} \coloneqq \min_{k \in \mc N} v^\star_k$ and $v^{\star}_{\max} \coloneqq \max_{k \in \mc N} v^\star_k$ are the smallest and largest magnitude set-points, and $\lambda_2(L)$ is the second smallest eigenvalue of the graph Laplacian $L$.
\end{condition}
If the graph $\mc G$ is connected (i.e.,  $\lambda_2(L)>0$), Condition \ref{cond.stab} can always be satisfied by a suitable choice of control gains and set-points. We refer to Proposition \ref{cond.stab.power} for a more detailed discussion.
We can now state the main result of the manuscript.
\begin{theorem}{\bf(Almost global stability of $\boldsymbol{\mc T}$)}\label{thm.main}
Consider set-points $p^\star_k$, $q^\star_k$, $v^\star_k$, steady-state angles $\theta^\star_{jk}$, a stability margin $c$, and control gains $\alpha$ and $\eta$, such that Condition~\ref{cond.stab} holds. Then, the dynamics \eqref{eq.closed.loop.rot} are almost globally asymptotically stable with respect to $\mathcal{T}$, and the origin $\mathbbl{0}_{n}$ is an exponentially unstable equilibrium.
\end{theorem}
Theorem~\ref{thm.main} guarantees almost global asymptotic stability of the set $\mathcal{T}$ of desired equilibria in the rotating frame (corresponding to a harmonic solution with the desired power flows).
In the literature on virtual oscillator control for power inverters~{\cite{johnson2014synchronization,johnson2016synthesizing,MS-FD-BJ-SD:14b,LABT-JPH-JM:13}} it is shown that, for connected graphs, global synchronization can be obtained for the trivial power flow solution with {identical angles and voltages}, i.e., $\theta_{jk}^\star=0$, $v^\star_k=v^\star$. Moreover, in the context of droop-controlled microgrids, \cite{vorobev2017high,vorobev2017cdc} provide stability conditions which account for the electromagnetic transients of the transmission lines and are valid locally around the trivial power flow solution. Condition~\ref{cond.stab} (together with Theorem~\ref{thm.main}) ensures that a larger set of power-flow solutions (with nonzero power flows) is almost globally stable and considers the destabilizing effects of the electromagnetic transients of the transmission lines. We remark that Condition \ref{cond.stab} is not overly conservative (see Section \ref{sec.threebusboundary}). Let us now provide a power-system interpretation of Condition~\ref{cond.stab}.
\begin{proposition}{\bf(Interpretation of the stability condition)}\label{cond.stab.power}
 Condition~\ref{cond.stab} is satisfied if the steady-state angles $\theta^\star_{jk}$, the set-points $p^\star_k$, $q^\star_k$, $v_k^\star$, and the steady-state branch powers $p^\star_{jk}$, $q^\star_{jk}$ satisfy Condition~\ref{cond.consistent}, $|\theta^\star_{jk}| \leq \tfrac{\pi}{2}$ holds for all $(j,k) \in \mc N \times \mc N$, and for all $k \in \mc N$, the stability margin $c \in \mathbb{R}_{>0}$, the network parameters $\norm{Y_{jk}}$, and the gains $\eta \in \mathbb{R}_{>0}$, $\alpha \in \mathbb{R}_{>0}$ satisfy
 \begin{align*}
  \sum_{j : (j,k) \in \mc E} \frac{\cos(\kappa)}{v^{\star2}_k} \left|p^\star_{jk}\right|\!+\!\frac{\sin(\kappa)}{v^{\star2}_k} \left| q^\star_{jk} \right|\!+\! \alpha \leq 
\frac{v^{\star2}_{\min}}{2 v^{\star2}_{\max}}\lambda_2( L)-c,\\
\eta < \frac{c}{2 \rho d_{\max} (c  + 5 \max_{k \in \mc N} s^\star_k  {v^\star_k}^{-2} + 10 d_{\max})},
 \end{align*}
where $s^\star_k \!\coloneqq\! \sqrt{p^{\star2}_k+q^{\star2}_k}$ is the apparent steady-state power injection, and $d_{\max} \!\coloneqq\! \max_{k \in \mc N} \sum\nolimits_{j:(j,k) \in \mc E} \norm{Y_{jk}}$ is the maximum weighted node degree of the transmission network graph.
\end{proposition}
{The proof is given in the Appendix. The stability conditions confirm and quantify well-known engineering insights. Broadly speaking, they require that  the network is not to heavily loaded and that there is a sufficient time-scale separation between the inverter dynamics and line dynamics.} The first inequality in Proposition~\ref{cond.stab.power} bounds the achievable steady-state power transfer in terms of the connectivity $\lambda_2(L)$ of the graph of the transmission network, the stability margin $c$ that accounts for the effect of the line dynamics, the control gain $\alpha$, and the voltage set-points $v^\star_{k}$ and shows that the achievable power transfer can be increased by increasing all voltage set-points $v^\star_{k}$, by reducing the gain $\alpha$ of the voltage regulator, or by upgrading transmission lines, i.e., increasing $\norm{Y_{jk}}$ and hence $\lambda_2(L)$. The second inequality in Proposition~\ref{cond.stab.power} provides a bound on the control gain $\eta$. This bound ensures that the line dynamics do not destabilize the system. We note, that for operating points that satisfy Proposition~\ref{cond.stab.power}, the term $\max_{k \in \mc N} s^\star_k {v^\star_k}^{-2}$ is much smaller than $d_{\max}$ and can be neglected for the purpose of this discussion. 

The controller \eqref{eq.control.law} achieves synchronization by inferring information about the voltage angle differences through the local measurements of the currents $\underline{i}_{o,k}$ (see Section \ref{sec.control.cl.design}). Thus, the time constant $\rho$ of the transmission lines can be interpreted as the propagation delay of the information on the phase angles. 
Loosely speaking, the controller cannot act faster than it can observe information through the network. Therefore, longer time-constants $\rho$ require a lower gain $\eta$. 

Moreover, the overall system gain is a combination of the gain of the transmission network, corresponding to the admittances $\norm{Y_{jk}}$ and the controller gains $\eta$. Upgrading or adding a line can increase $d_{\max}$ more than $\lambda_2(L)$. Thus, to keep the closed-loop gain constant, the controller gain $\eta$ needs to decrease. In other words, {the sufficient stability conditions of Theorem \ref{thm.main} and Proposition \ref{cond.stab.power} indicate that the system may become unstable if transmission lines are added, shortened, or upgraded. This observation is verified in Section \ref{sec:three:adm} using a numerical example.} The same effect was observed in \cite{vorobev2017cdc} for droop controlled microgrids. In the remainder of this section we prove Theorem~\ref{thm.main}.

\subsection{Singular perturbation theory}\label{subsec.sing.pert}
In the following we apply tools from singular perturbation theory (see \cite[{Ch. 11.5}]{K02}) to explicitly construct a Lyapunov function that establishes convergence of the dynamics~\eqref{eq.closed.loop.rot} to the set $\mc T_0$ and allows us to show that the origin $\mathbbl{0}_{n}$ is an unstable equilibrium. By replacing the dynamics \eqref{eq.lines.cl.rot} of the transmission lines with its steady-state map $i^s(v)$, we obtain the \emph{reduced-order system}
\begin{align}\label{eq.voltage.control.kh}
\begin{split}
\ddt v = f_v(v,i^s(v)) = {\eta}\left[( \mc K - \mc L)v  + \alpha \Phi(v){v}\right ] ,
\end{split}
\end{align}
which describes the voltage dynamics under the assumption that the line currents are at their quasi-steady-state, i.e., $i_o = i^s_o(v)$. This results in $ \mc R(\kappa) i_o = \mc R(\kappa) \mc Y v = \mc L v$. Moreover, note that $(\mc K - \mc L) v = e_\theta(v)$ which implies $(\mc K - \mc L) v = 0$ for all $v \in \mc S$. Finally, we denote the difference between the line currents and their steady-state value as $y = i - i^s(v)$ and define the \emph{boundary-layer system}
\begin{align}
\ddt y \bigg\vert_{\ddt v = \mathbbl{0}_{2N}} = f_i(v,y+i^s(v)) = {-L^{-1}_T Z_T y} \label{eq.current.control.kh}
\end{align}
where $v$ is treated as a constant. In the remainder, we follow a similar approach to~\cite[Sec. 11.5]{K02} and obtain stability conditions for the full system~\eqref{eq.closed.loop.rot} based on the properties of~\eqref{eq.voltage.control.kh} and~\eqref{eq.current.control.kh}. In contrast to~\cite{K02} we address stability with respect to the set $\mc T_0$ and not a single equilibrium. 

\subsection{Lyapunov function for the reduced-order system}\label{sec.power.systems.steady}
In this section, we provide a Lyapunov function for the reduced-order system~\eqref{eq.voltage.control.kh}. Given the voltage set-points $v^\star_k$ and relative steady-state angles $\theta^\star_{k1}$ for all $k \in \mc N$, we define the matrix $S\coloneqq[v_1^\star R(\theta_{11}^\star)^\mathsf{T} \hdots v_N^\star R(\theta_{1N}^\star)^\mathsf{T}]^\mathsf{T}$ and the projector $P_S \coloneqq (I_{2N} - \frac{1}{\sum v_i^{\star 2}}S S^\mathsf{T})$ onto the nullspace of ${\mc S}$. We now define the Lyapunov function candidate $V:\R^n\to\R_{\ge0}$ for the reduced-order system as
\begin{align}\label{eq.v}
V(v) \coloneqq \frac{1}{2} v^\mathsf{T}P_Sv +  {\frac{1}{2}}\eta\alpha\alpha_1\sum_{k=1}^N \left( \frac{{v_k^\star}^2 - \norm{v_k}^2}{v_k^\star}\right)^2,
\end{align}
where $\alpha \in \mathbb{R}_{>0}$ is the voltage controller gain and, given $c \in \mathbb{R}_{>0}$ such that Condition \ref{cond.stab} holds, the constant $\alpha_1 $ is given by
\begin{align} \label{thm1.strictbound}
\alpha_1 \coloneqq \frac{ c  }{ 5\eta \norm{\mc K -\mc L}^2}.
\end{align}
Moreover, we define the function $\psi: \R^{2N} \to \R_{\geq 0}$ as
\begin{align}\label{eq.voltage.psi.def}
\psi(v) &\coloneqq \eta  \left(\|\mc K-\mc L\|  \|v\|_{\mc S}  + \alpha  \|\Phi(v)v\| \right).
\end{align}
We require the following preliminary results
\begin{lemma} \label{lem.positivsatz.1}
The following inequality holds for all $v\in\R^{2N}$
\begin{align}\label{eq.claim1}
v^\mathsf{T} P_S \Phi(v) v \leq v^\mathsf{T} P_S v = \|v\|_{\mc S}^2.
\end{align}
\end{lemma}
\begin{lemma} \label{lem.redlfdecr}
Consider steady-state angles $\theta^\star_{jk}$, $\alpha$, and $c$ such that Condition~\ref{cond.stab} is satisfied. For all $v\in\R^{2N}$, it holds that
 \begin{align}\label{eq.ass.reformulation}
   v^\mathsf{T} P_S ( \mc K - \mc L + \alpha I_{2N}) v \leq - c \,\|v\|_{\mc S}^2.
 \end{align}
\end{lemma}
The proofs are given in the appendix and use the same technique as in~\cite[Prop. 6, Prop. 7]{colombino2017global2}. In the following {proposition}, we show {that} the function $V$ is a Lyapunov function for the reduced-order system~\eqref{eq.voltage.control.kh}.
\begin{proposition}{\bf(Lyapunov function for the reduced-order system)}\label{thm.lyap.red.system}
Consider $V(v)$ defined in \eqref{eq.v} and set-points $p^\star_k$, $q^\star_k$, $v^\star_k$, steady-state angles $\theta^\star_{jk}$, $\alpha$ and $c$ such that Condition~\ref{cond.stab} holds. For any $\eta \in \mathbb{R}_{>0}$, there exists $\chi^{\scriptscriptstyle{V}}_1, \chi^{\scriptscriptstyle{V}}_2 \in \mathscr{K}_{\infty}$ such that
\begin{align}\label{eq.Vbound}
\chi^{\scriptscriptstyle{V}}_1(\norm{v}_{\mc S \cap \mc A}) \leq V(v) \leq  \chi^{\scriptscriptstyle{V}}_2(\norm{v}_{\mc S \cap \mc A})
\end{align}
holds for all $v \in \R^{2N}$. Moreover, for all $v \in \R^{2N}$ the derivative of $V$ along the trajectories of the reduced-order system~\eqref{eq.voltage.control.kh} satisfies
\begin{align} \label{eq.voltage.ss.Vdot.ideal}
\ddt V \coloneqq \frac{\partial V}{\partial v}  f_v(v,i^s(v)) &\leq -\alpha_1 \psi(v)^2.
\end{align}
\end{proposition}
\begin{IEEEproof}
By construction the function $V(v)$ is positive definite with respect to $\mc S \cap \mc A$, i.e., $V=0$ for all $v \in \mc S \cap \mc A$ and $V(v)>0$ otherwise, and radially unbounded with respect to $\mc S \cap \mc A$, i.e., $V(v) \to \infty$ for $\norm{v}_{\mathcal{A}} \to \infty$ and $V(v) \to \infty$ for $\norm{v}_{\mathcal{S}} \to \infty$. Moreover $\mc S \cap \mc A$ is a compact set. Using these properties and replacing $\norm{\cdot}$ with the point to set distance $\norm{\cdot}_{\mc S \cap \mc A}$ in the steps outlined in \cite[p. 98]{H67}, it directly follows that there exist $\mathscr{K}_\infty$-functions $\chi^{\scriptscriptstyle{V}}_1$ and $\chi^{\scriptscriptstyle{V}}_2$ such that \eqref{eq.Vbound} holds for all $v \in \R^{2N}$.
Next, we can write the derivative of $V(\cdot)$ along the trajectories of~\eqref{eq.voltage.control.kh} as
\begin{align}\label{eq.lyap.decrease.one}
\begin{split}
\ddt V & = \eta v^\mathsf{T} P_S \left( ( \mc K - \mc L)v  + \alpha \Phi(v)v\right)\\
& - 2 \eta^2 \alpha\,\alpha_1\,v^\mathsf{T} \Phi(v)  \left(( \mc K - \mc L)v  + \alpha \Phi(v)v\right).
\end{split}
\end{align}
Using Lemma~\ref{lem.positivsatz.1}, we can bound~\eqref{eq.lyap.decrease.one} as
\begin{align}\label{eq.lyap.decrease.lemm1}
\begin{split}
\ddt V \le &\eta v^\mathsf{T} P_S ( \mc K - \mc L + \alpha I_{2N}) v \\
& - 2 \eta^2 \alpha\,\alpha_1\,v^\mathsf{T} \Phi(v)  \left(( \mc K - \mc L)v  + \alpha \Phi(v)v\right).
\end{split}
\end{align}
Using Lemma~\ref{lem.redlfdecr} we obtain
\begin{align*}\label{eq.lyap.decrease.two}
\ddt V  \le -\eta c \norm{v}^2_{\mc S} - 2\eta^2 \alpha \alpha_1 v^{\mathsf{T}} \Phi(v)\!  \left(( \mc K - \mc L)v  + \alpha \Phi(v)v\right).
\end{align*}
Because $(\mc K - \mc L) v =0$ for all $v\in\mc S$ and $P_S$ is the projector onto $\mc S^\perp$ it holds that $\mc K-\mc L = (\mc K-\mc L)P_{S}$, and we obtain
\begin{align*}\label{eq.lyap.decrease.three}
\ddt V \le & -\eta c\, \|v\|_{\mc S}^2 -2 \eta^2 \alpha^2 \alpha_1 \| \Phi(v)v\|^2\\
& + 2 \eta^2 \alpha \alpha_1\norm{\Phi(v) v} \norm{\mc K - \mc L} \norm{v}_S.
\end{align*}
Next, to show that \eqref{eq.voltage.ss.Vdot.ideal} holds, we show that for $\alpha_1$ according to \eqref{thm1.strictbound} and for all $v \in \mathbb{R}^{2N}$, the inequality
\begin{align}\label{eq.inequality}
\begin{split}
 &-\eta c\, \|v\|_{\mc S}^2 - 2\eta ^2\alpha^2 \alpha_1 \| \Phi(v)v\|^2 \\
 &+ 2\, \eta^2 \alpha \alpha_1\norm{\Phi(v) v} \norm{\mc K - \mc L} \norm{v}_S\\
 &\le - \alpha_1 \eta^2(\|\mc K-\mc L\| \|v\|_\mc S +\alpha\| \Phi(v) v\| )^2
 \end{split}
\end{align}
holds. By matching the coefficients of the r.h.s and the l.h.s., it can be seen that~\eqref{eq.inequality} holds for all $v \in \mathbb{R}^{2N}$ if %
\begin{align}\label{eq.condition}
\begin{bmatrix}
\|v\|_{\mc S}\\
\| \Phi(v)v\|
\end{bmatrix}^\mathsf{T}
Q
\begin{bmatrix}
\|v\|_{\mc S}\\
\| \Phi(v)v\|
\end{bmatrix}\ge0,\quad \forall v\in \R^{2N},
\end{align}
holds for the matrix $Q$ defined as
\begin{align}\label{eq.cond.Q}
Q\coloneqq  \begin{bmatrix}
\eta c - \alpha_1\eta^2\|\mc K-\mc L\|^2 & -2\alpha_1\alpha\eta^2\|\mc K-\mc L\|\\
\star & \eta^2\alpha^2\alpha_1
\end{bmatrix}. 
\end{align}
Using the Schur complement and $\alpha_1 >0$ it follows that $Q\succ 0$ if $\eta c - 5\alpha_1\eta^2\|\mc K-\mc L\|^2 \ge0$. Thus, \eqref{eq.condition} is satisfied for $\alpha_1$ defined in \eqref{thm1.strictbound} and the proposition directly follows.
\end{IEEEproof}
\subsection{Lyapunov function for the Boundary layer system}\label{sec.power.systems.boundary}
In this section, we provide a Lyapunov function for the boundary layer system~\eqref{eq.current.control.kh}. To this end, we define $y_o \coloneqq i_o - i^s_o(v) = \mc B y$, the matrix $B_n \in \R^{M \times M_0}$ whose columns span the nullspace of $B$, $\mc B_n \coloneqq B_n \otimes I_2$, and $y_n = \mc B^\mathsf{T}_n L_T y$. The fact that $\mc G$ is connected implies $\rank(B)=N-1$ and it follows from the rank-nullity theorem that $M_0 = M-N+1$. We define the Lyapunov function $W: \R^{2M} \to \R_{\geq 0}$
\begin{align} \label{eq.w}
W(y) \coloneqq \frac{\rho}{2} (y^\mathsf{T} \mc B^\mathsf{T} \mc B y + y^\mathsf{T} L_T \mc B_n \mc B^\mathsf{T}_n L_T y).
\end{align}
\begin{proposition}{\bf(Lyapunov function for the boundary layer system)} \label{prop.deriv.W}
There exists $\chi^{\scriptscriptstyle{W}}_1, \chi^{\scriptscriptstyle{W}}_2 \in \mathscr{K}_{\infty}$ such that $ \chi^{\scriptscriptstyle{W}}_1(\norm{y}) \leq W(y) \leq \chi^{\scriptscriptstyle{W}}_2(\norm{y})$ holds for all $y \in \R^{2M}$. Moreover, for every $v \in \R^{2N}$, and all $y \in \R^{2M}$, the derivative of $W(y)$ along the trajectories of~\eqref{eq.current.control.kh} satisfies
\begin{align*}
\ddt W  \big\vert_{\ddt v = \mathbbl{0}_{2N}} = \frac{\partial W}{\partial y} f_i(v,y+i^s(v)) \le -\norm{y_o}^2 - \norm{y_n}^2.
\end{align*}
\end{proposition}
\begin{IEEEproof}
Because $\mc B_n$ has full rank, it holds that $\mc B^\mathsf{T}_n L_T \mc B_n \succ 0$. Next, we can parametrize all $y$ such that $\mc B y \!=\!\mathbbl{0}_{2N}$ using $\xi \in \R^{2M_0}$ and $y = \mc B_n \xi$, and it follows that $\xi^\mathsf{T} (\mc B^\mathsf{T}_n L_T \mc B_n)^\mathsf{T}(\mc B^\mathsf{T}_n L_T \mc B_n) \xi >0$ for all $\xi$.
 Thus, $W(y)$ is positive definite and the first statement directly follows. Using $L^{-1}_T Z_T = \tfrac{1}{\rho}I_{2M}+\omega_0 \mc J_M$, one obtains
\begin{align*}
\ddt W \big\vert_{\ddt v = \mathbbl{0}_{2N}} =& -y^\mathsf{T}_o y^{\phantom{\mathsf{T}}}_o - \rho y^\mathsf{T} \mc B^\mathsf{T} \mc B   \omega_0 \mc J_M y\\
&- y^\mathsf{T}_n y^{\phantom{\mathsf{T}}}_n - \rho y^\mathsf{T} L_T \mc B_n \mc B^\mathsf{T}_n L_T \omega_0 \mc J_M y.
\end{align*}
Using the mixed-product property of the Kronecker product and $J\!=\!-J^\mathsf{T}\!$, it can be shown that $\mc B^\mathsf{T}\! \mc B  \mc J_M \!=\! ((B^\mathsf{T}\! B) \otimes I_2)(I_M \otimes J) \!=\! (B^\mathsf{T}\! B) \otimes J \!=\! -((B^\mathsf{T}\! B) \otimes J)^\mathsf{T}\!$ is skew-symmetric. The same approach can be used to show that $L_T \mc B_n \mc B^\mathsf{T}_n L_T \mc J_M$ is skew-symmetric and the proposition follows.
\end{IEEEproof}

\subsection{Proof of the main result}
Proposition~\ref{thm.lyap.red.system} establishes that the Lyapunov functions $V(v)$ of the reduced-order system~\eqref{eq.voltage.control.kh}, which describes the voltage dynamics assuming the transmission line currents are in steady state, is decreasing. Moreover, Proposition~\ref{prop.deriv.W} shows that the Lyapunov function $W(y)$ of the boundary layer system~\eqref{eq.current.control.kh}, which describes the difference between the actual transmission line currents and their steady state for constant voltages (in the rotating reference frame), is decreasing. 

In this section, we use these results to construct a Lyapunov function candidate for the overall system~\eqref{eq.closed.loop.rot} as a convex combination of the functions $V(v)$ and $W(y)$ introduced in~\eqref{eq.v} and~\eqref{eq.w}. We recall that $x \coloneqq (v,i) \in \R^n$, introduce a constant $d \in \mathbb{R}_{(0,1)}$, and define the Lyapunov function candidate $\nu: \R^n \to \R_{\geq0}$ for the full-order system~\eqref{eq.closed.loop.rot} as
\begin{align}\label{eq.lyap.full}
 \nu(x)\coloneqq d W (i-i^s(v))+ (1-d)V(v).
\end{align}
The following propositions bound the terms which arise from taking the time-derivative of $\nu(x)$ for the full-order system~\eqref{eq.closed.loop.rot} instead of the reduced-order system~\eqref{eq.voltage.control.kh} and boundary layer system~\eqref{eq.current.control.kh}. We will use these bounds to define a constant $d \in \mathbb{R}_{(0,1)}$ that ensures that $\nu(x)$ is decreasing.
\begin{proposition} \label{prop.beta1}
Let $\beta_1  \coloneqq \norm{\mc K -\mc L}^{-1}$. For all $v\in\R^{2N}$ and all $y \in \R^{2M}$ it holds that
\begin{align*}
\frac{\partial V}{\partial v}[f_v(v,y+i^s(v))-f_v(v,i^s(v))]\le  \beta_1\psi(v)\|y_o\|. 
\end{align*}
\end{proposition}
\begin{proposition} \label{prop.beta2}
Consider $\beta_2 \coloneqq \rho \norm{\mc Y}$ and $\gamma \coloneqq \eta \rho \norm{\mc Y}$. For all $v\in\R^{2N}$ and $y \in \R^{2M}$ it holds that
\begin{align*}
-\frac{\partial W}{\partial y}\frac{\partial i^s}{\partial v} f_v(v,y+i^s(v))\leq
 \beta_2 \psi(v) \|y_o\| + \gamma \|y_o\|^2 .
\end{align*}
\end{proposition}
The proofs are provided in the Appendix.
The next result bounds the decrease of $\nu$ along the trajectories of the \emph{full-order system}~\eqref{eq.closed.loop.rot}.
\begin{proposition}\label{propo.nudecr}
Let $d \coloneqq \frac{\beta_1}{\beta_1+\beta_2} \in \mathbb{R}_{(0,1)}$. Under Condition~\ref{cond.stab}, there exist a function $\chi_3\in \mathscr{K}$ such that 
\begin{align}\label{eq.bound.decrease.nu}
\ddt \nu \coloneqq \frac{\partial \nu}{\partial v} f_v(v,i) + \frac{\partial \nu}{\partial v} f_i(v,i) \le - \chi_3(\|x\|_{\mc T_0}).
\end{align}
\end{proposition}
\begin{IEEEproof}
It follows from Propositions~\ref{thm.lyap.red.system}~-~\ref{prop.beta2} 
\begin{align*}
&\ddt \nu = d \bigg(\! \frac{\partial W}{\partial y} f_i(v,y+i^s(v)) - 
\frac{\partial W}{\partial y}\frac{\partial i^s}{\partial v} f_v(v,y+i^s(v))\!\!\bigg)+\\
&(1\!-\!d)\bigg(\!\frac{\partial V}{\partial v} \!f_v(v,i^s(v)) \!+\! \frac{\partial V}{\partial v}[f_v(v,y\!+\!i^s(v))\!-\!f_v(v,i^s(v))]\!\bigg)\!,
\end{align*}
$y_o = i_o-i^s_o(v)$, and $y_n = \mc B^\mathsf{T}_n L_T y$, that
\begin{align*}
&\ddt \nu \leq - d \| y_o \|^2 - d \| y_n \|^2  + d \beta_2 \psi(v) \| y_o\|  + d\gamma \|  y_o\|^2\\ 
&\hphantom{\ddt \nu \leq} -(1-d)\alpha_1\psi(v)^2 +(1-d)\beta_1\psi(v) \| y_o\|   \\
\leq & -\!\begin{bmatrix}
\psi(v)\\
\| y_o\|\\
\| y_n \|
\end{bmatrix}^{\!\mathsf{T}}\!\!
\underbrace{\begin{bmatrix}
 (1-d)\alpha_1\!\! & -\frac{(1-d)\beta_1+ d\beta_2}{2} & 0\\
 \star& (1-\gamma) d & 0 \\
 0 & 0 & {d}
\end{bmatrix}}_{\eqqcolon \mc M}\!
 \begin{bmatrix}
\psi(v)\\
\| y_o\|\\
\| y_n \|
\end{bmatrix}
\end{align*}
Following similar steps to~\cite[p. 452]{K02}, we let $d = \frac{\beta_1}{\beta_1+\beta_2}$ and obtain the following condition for $\mc M$ to be positive definite
\begin{align}\label{eq.singper.cond}
1  < \frac{{\alpha_1} }{\alpha_1 \gamma + \beta_1 \beta_2}\implies \mc M\succ 0.
\end{align}
Using the expressions for $\alpha_1,\beta_1,\beta_2$ and $\gamma$ derived in Propositions~\ref{thm.lyap.red.system}~-~\ref{prop.beta2}, condition~\eqref{eq.singper.cond} is satisfied if
\begin{align} \label{eq.eta.def}
\eta < \frac{c }{\rho \norm{\mc Y} (c  + 5 \|\mc K -\mc L\|)}.
\end{align}
Therefore $\ddt \nu(x)<0$ for all $(v,i)\notin \mc T_0$. Moreover, because ${(\psi(v)^2+\|y_o\|^2+\|y_n\|^2)}$ is both positive definite and radially unbounded with respect to $\mc T_0$, and $\mc T_0$ is compact, the steps in~\cite[p. 98]{H67} can be used to show that there exists a function $\chi_3\in \mathscr{K}$ such that~\eqref{eq.bound.decrease.nu} holds.
\end{IEEEproof}

We are now ready to prove Theorem~\ref{thm.main}.\\
\emph{Proof of Theorem~\ref{thm.main}:}
We apply Theorem \ref{thm.AGAS} with $\mc C = \mc T$ and $\mc U = \{\mathbbl{0}_n\}$. According to Proposition~\ref{thm.lyap.red.system} and \ref{prop.deriv.W} there exists $\mathscr{K}_{\infty}$-functions $\chi^{\scriptscriptstyle{V}}_1$, $\chi^{\scriptscriptstyle{V}}_2$, $\chi^{\scriptscriptstyle{W}}_1$, and $\chi^{\scriptscriptstyle{W}}_2$, such that { the functions $\tilde\chi_i = (1-d) \chi^{\scriptscriptstyle{V}}_i + d \chi^{\scriptscriptstyle{W}}_i$, $i\in\{1,2\}$ are positive definite and radially unbounded w.r.t the compact set $\mc T$ and $\tilde\chi_1(\norm{v}_{\mc S\cap\mc A },\norm{y})\le \nu(x)\le \tilde\chi_2(\norm{v}_{\mc S\cap\mc A },\norm{y})$, therefore, from~\cite[p. 98]{H67} there exist $\chi_1,\chi_2\in\mathscr{K}_{\infty}$ such that}
$\chi_1(\norm{x}_{\mc T})\le \nu(x)\le \chi_2(\norm{x}_{\mc T})$ holds for all $x\in \mathbb{R}^{n}$. Moreover, Proposition~\ref{propo.nudecr} guarantees the existence of a function $\chi_3 \in \mathscr{K}$ such that
$\ddt \nu \le - \chi_3(\|x\|_{\mc T_0})$ holds for all $x\in \mathbb{R}^{n}$. It remains to show that the region of attraction $\mc Z_{\mc U} = \mc Z_{\{\mathbbl{0}_n\}}$ has zero Lebesgue measure. 
{To this end, consider the linearized dynamics $\ddt x \!=\! A_{\mathbbl{0}} x$ with $A_{\mathbbl{0}} \coloneqq \tfrac{\partial}{\partial x} f(x) \vert_{x = \mathbbl{0}_{n}}$ and the reduced-order linearized dynamics $\ddt v = A_{\mathbbl{0},v} v$ with $A_{\mathbbl{0},v} \coloneqq \tfrac{\partial}{\partial v}f_v(v,i^s(v))\vert_{v = \mathbbl{0}_{2N}}$. Moreover, we define $\psi_{\mathbbl{0}}(v) \coloneqq  \eta ( \norm{\mc K - \mc L} \norm{v}_{\mc S} + \alpha \norm{v})$ and the quadratic approximation of the previous Lyapunov function $\nu_{\mathbbl{0}}(x) \coloneqq \frac{1}{2}x^\mathsf{T}  \tfrac{\partial^2}{\partial^2 x} \nu(x)\vert_{x=\mathbbl{0}_n} x$. Following analogous steps as in the proofs of Propositions \ref{thm.lyap.red.system} and \ref{propo.nudecr}, it can be shown that $\frac{\partial V_{\mathbbl{0}}}{\partial v} A_{\mathbbl{0},v} v \leq -\alpha_1 \psi_{\mathbbl{0}}(v)$ holds, and that there exists positive constants $d \in \mathbb{R}_{(0,1)}$ and $\chi_{\mathbbl{0}} \in \mathbb{R}_{>0}$ such that $\ddt \nu_{\mathbbl{0}} \coloneqq \frac{\partial \nu_{\mathbbl{0}}}{\partial x} A_{\mathbbl{0}} \leq - \chi_{\mathbbl{0}} (\psi_{\mathbbl{0}}^2 + \norm{y_o}^2 + \norm{y_n}^2)$ holds.  Next, we show that $A_{\mathbbl{0}}$ and cannot have eigenvalues with zero real part. If $A_{\mathbbl{0}}$ has eigenvalues with zero real part there exists $x_0 \neq \mathbbl{0}_n$ such that the solution $\varphi_{A_{\mathbbl{0}}}(t,x_0)$ of $\ddt x_\delta = A_{\mathbbl{0}} x_\delta$ 
remains bounded for all $t \in \mathbb{R}_{\geq 0}$ and does not converge to the origin. However, $\nu_{\mathbbl{0}}(\varphi_{A_{\mathbbl{0}}}(t,x_0))$ is strictly decreasing in $t$ for all $x_0 \neq \mathbbl{0}_n$ and $\nu_{\mathbbl{0}}$ is quadratic and not bounded from below. Thus, it either holds that $\norm{\varphi_{A_{\mathbbl{0}}}(t,x_0)} \to 0$ as $t \to \infty$ or $\nu_{\mathbbl{0}}(\varphi_{A_{\mathbbl{0}}}(t,x_0)) \to -\infty$ and $\norm{\varphi_{A_{\mathbbl{0}}}(t,x_0)} \to -\infty$ as $t \to \infty$, i.e., there exists no initial condition $x_0 \neq \mathbbl{0}_n$ for which $\varphi_{A_{\mathbbl{0}}}(t,x_0)$ remains bounded and does not converge to the origin. Therefore, $A_{\mathbbl{0}}$ cannot have eigenvalues with zero real part. Moreover, for all $x_0=(v_0,i^s(v_0))$ with $v_0 \in \mc S \setminus \{\mathbbl{0}_{2N}\}$ it holds that $\nu_{\mathbbl{0}}(\mathbbl{0}_n) > \nu_{\mathbbl{0}}(x_0)$, and it follows that $\nu_{\mathbbl{0}}(\varphi_{A_{\mathbbl{0}}}(t,x_0)) \to -\infty$ and $\norm{\varphi_{A_{\mathbbl{0}}}(t,x_0)} \to -\infty$ as $t \to \infty$, i.e., the origin is an unstable equilibrium of $\ddt x = A_{\mathbbl{0}} x$,  and at least one eigenvalue of $A_{\mathbbl{0}_n}$ has positive real part. Because the right-hand side of \eqref{eq.closed.loop.rot} is continuously differentiable, the region of attraction $\mc Z_{\mc U} = \mc Z_{\{\mathbbl{0}_n\}}$ of the origin has zero Lebesgue measure (see \cite[Prop. 11]{Monzon2006}). Moreover, it holds that $\mathbbl{0}_n \notin \mc C$ and the theorem directly follows from Theorem~\ref{thm.AGAS}.
\hfill\IEEEQED }

\section{Power systems implications and test-cases}\label{sec.example}
In this section, two power systems test-cases are used to illustrate the results. We use a three-bus system to {illustrate the impact of the line parameters and} to illustrate and discuss the stability boundaries obtained from Theorem \ref{thm.main}. Then, we consider a more realistic test-case based on the IEEE 9-bus system, and we show that the system behaves as expected even when some of the technical assumptions do not hold.

{
\subsection{The impact of line admittances on stability}\label{sec:three:adm}
We now investigate the role line admittances play in the stability of an inverter-based transmission system using the system depicted in Figure~\ref{fig:threebus} consisting of three inverters and three transmission lines. The base power is $1$ GW, the base voltage $320$ kV, the transmission line connecting inverter 1 to inverter 2 is $125$ km long, and the transmission lines connecting inverter 1 and inverter 3, as well as inverter 2 and inverter 3, are $25$ km long. The line resistance is $0.03$ Ohm/{km} and the line reactance is $0.3$ Ohm/km (at $\omega_0 = 50$ Hz), i.e., the $\ell/r$ ratio of the transmission lines is $\omega_0\rho=10$. We use dVOC to control the inverters with the set-points (in p.u.), $v_k^\star=1$, $p^\star_1=-0.52$, $p^\star_2= -0.19$, $p^\star_3= 0.71$, $q^\star_1=0.06$, $q^\star_2=0.021$, and $q^\star_3=-0.06$, $\eta = 3\cdot 10^{-3}$ and $\alpha=5$.

Theorem~\ref{thm.main} and Proposition~\ref{cond.stab.power} indicate that the system may become unstable if the admittance of individual lines is increased by, e.g., adding or upgrading transmission lines. To validate this insight, we  vary the admittance of individual lines (keeping $\ell/r$ constant), recompute the steady-state given by $\theta^\star_{jk}$ and $v^\star_k=1$ that corresponds to the power injections specified above, linearize the system at this steady-state, and compute the minimum damping ratio $\zeta_{\min}$, defined by 
\begin{align*}
 \zeta_{\min} \coloneqq \min_{k \in \{2,\ldots,n\}} \frac{\operatorname{Re}(\lambda_k)}{\sqrt{\operatorname{Re}(\lambda_k)^2 + \operatorname{Im}(\lambda_k)^2}},
\end{align*}
where $\operatorname{Re}(\lambda_k)$ and $\operatorname{Im}(\lambda_k)$ denote the real part of the $k$-th eigenvalue of the linearized system and we exclude the eigenvalue $\lambda_1=0$ that corresponds to the rotational invariance of the system. A larger damping ratio corresponds to a well damped system, and the damping ratio is negative if the system is unstable.
\begin{figure}[b!!!]
\definecolor{mycolor1}{rgb}{0.00000,0.44700,0.74100}%
\definecolor{mycolor2}{rgb}{0.85000,0.32500,0.09800}%
\definecolor{mycolor3}{rgb}{0.92900,0.69400,0.12500}%
\newlength\figH
\newlength\figW
\setlength{\figH}{1.6cm}
\setlength{\figW}{7.5cm}
\begin{center}
\begin{circuitikz}[american voltages]

\ctikzset{bipoles/resistor/height=0.15}
\ctikzset{bipoles/resistor/width=0.4}

\ctikzset{bipoles/generic/height=0.15}
\ctikzset{bipoles/generic/width=0.4}

\ctikzset{bipoles/length=.6cm}

\coordinate (c) at (0.5,\figH+2cm);

\coordinate (I1) at   ($ (c) +  (0,0) $);
\coordinate (I2) at   ($ (c) + (6,0.5) $);
\coordinate (I3) at   ($ (c) + (2,-0.5) $);

\draw ($ (I1) + (0,-0.8) $) node [ground] (g1) {};
\draw (I1) to[sV, l=$v_2$,*-] (g1);

\draw ($ (I2) + (0,-0.8) $) node [ground] (g2) {};
\draw (g2) to[sV, l=$v_1$,-*] (I2);

\draw ($ (I3) + (0,-0.8) $) node [ground] (g3) {};
\draw (g3) to[sV, l=$v_3$,-*] (I3);

\node (I12) at ($(I1)!0.5!(I2)$) [label={[xshift=0cm, yshift=0cm,color=mycolor1]$125~\text{Km}$},color=mycolor1]{};
\draw 		(I1) 
to[R,color=mycolor1] (I12)
to [L,color=mycolor1] (I2);

\node (I13) at ($(I1)!0.5!(I3)$)[label={[xshift=1.2cm, yshift=-0.2cm,color=mycolor2]$25~\text{Km}$},color=mycolor2]{};
\draw 		(I1) 
to[R,color=mycolor2] (I13)
to [L,-,color=mycolor2] (I3);

\node (I32) at ($(I3)!0.5!(I2)  $) [label={[xshift=0.3cm, yshift=-0.6cm,color=mycolor3]$25~\text{Km}$},color=mycolor3]{};
\draw 		(I3) 
to[R,color=mycolor3] (I32)
to [L,color=mycolor3] (I2);
\end{circuitikz}%
\vspace{-1em}
\end{center}
\caption{{Three-bus transmission system}. \label{fig:threebus}}
\vspace{1em}
%
%
\definecolor{mycolor1}{rgb}{0.00000,0.44700,0.74100}%
\definecolor{mycolor2}{rgb}{0.85000,0.32500,0.09800}%
\definecolor{mycolor3}{rgb}{0.92900,0.69400,0.12500}%
\begin{tikzpicture}

\begin{axis}[%
width=0.951\figW,
height=\figH,
at={(0\figW,0\figH)},
scale only axis,
xmin=0.027,
xmax=0.4,
xlabel style={font=\color{white!15!black}},
xlabel={{$\norm{Y_l}$ [S]}},
ymin=-0.002,
ymax=0.104,
ylabel style={font=\color{white!15!black}},
ylabel={{$\zeta_{\min}$}},
yticklabel style={
        /pgf/number format/fixed,
        /pgf/number format/precision=4
},
scaled x ticks=false,
xticklabel style={
        /pgf/number format/fixed,
        /pgf/number format/precision=4
},
scaled x ticks=false,
axis background/.style={fill=white},
xmajorgrids,
ymajorgrids,
legend style={legend cell align=left, align=left, draw=white!15!black},
legend pos={north east}
]
\addplot [color=black, dashed, line width=2.0pt, forget plot]
  table[row sep=crcr]{%
0.027	0.0686905244753528\\
0.3	0.0686905244753528\\
};
\addplot [color=mycolor1, line width=2.0pt]
  table[row sep=crcr]{%
0.266483545206184	-0.0126\\
0.245964249265391	-0.000643730619103833\\
0.226456739840894	0.0109111389236742\\
0.209816148734375	0.0208993011876092\\
0.195453733791248	0.0295980054894927\\
0.171917420298061	0.0440775064317775\\
0.153440314378176	0.0555796143178552\\
0.138549482181138	0.0649389062494727\\
0.132137735520844	0.0686916718480178\\
0.0463461217740715	0.0686892759389662\\
0.0267832196386049	0.0686905239493841\\
};
\addlegendentry{$\norm{Y_{12}}$}

\addplot [color=mycolor2, line width=2.0pt]
  table[row sep=crcr]{%
0.268632629863854	-0.00454931624232907\\
0.225726029260599	0.0193501920200473\\
0.194638013254649	0.0365335378805511\\
0.171076569544875	0.0492490401216558\\
0.15260351273956	0.0588646521923228\\
0.137731136497993	0.066263002592956\\
0.125500186152611	0.0720573248453967\\
0.115264355367115	0.0766126155942465\\
0.106572289224676	0.0802920131928096\\
0.0990992323583119	0.0832647206527771\\
0.0926055504658869	0.0856935741082887\\
0.0869105567206586	0.0877208014110448\\
0.0818754363061117	0.0894223361629954\\
0.077391781460777	0.0908643994315483\\
0.0733736979989307	0.0921059086951861\\
0.0664714687393176	0.0940887835942307\\
0.0607561648850959	0.0956161171776892\\
0.0559458661162243	0.0968189967478675\\
0.0518413847105683	0.0977875759351262\\
0.0467019370883999	0.0989281356887873\\
0.0424896055078776	0.09980460593372\\
0.0369369866062799	0.10088701878195\\
0.0319298116046428	0.10179314812356\\
0.0265343250722664	0.102702377963286\\
};
\addlegendentry{$\norm{Y_{23}}$}

\addplot [color=mycolor3, dashed, line width=2.0pt]
  table[row sep=crcr]{%
0.268632629863854	-0.00455379239129849\\
0.225726029260599	0.0193460780613097\\
0.194638013254649	0.0365302169996682\\
0.171076569544875	0.0492489957358219\\
0.15260351273956	0.0588633119094786\\
0.137731136497993	0.0662626394328985\\
0.125500186152611	0.072040756621817\\
0.115264355367115	0.0766139589203857\\
0.106572289224676	0.080279813079886\\
0.0990992323583119	0.0832544573250917\\
0.0926055504658869	0.0856955460227016\\
0.0869105567206586	0.0877314185279848\\
0.0818754363061117	0.0894263311159992\\
0.077391781460777	0.0908687154049502\\
0.0733736979989307	0.0921027059045584\\
0.0664714687393176	0.0940938117722183\\
0.0607561648850959	0.0956174406432777\\
0.0559458661162243	0.0968200522594909\\
0.050006997251579	0.0982096124513199\\
0.0452079947336945	0.099251179536554\\
0.0379282942981637	0.100706339038749\\
0.0326678876517852	0.101670116554997\\
0.0265343250722664	0.10270786197811\\
};
\addlegendentry{$\norm{Y_{13}}$}
\end{axis}
\end{tikzpicture}%
\vspace{-2em}
\caption{{Damping and stability of a three-bus transmission system as a function of the line admittances. The black dashed line indicates the minimum damping ratio for the system in the original configuration described in Section \ref{sec.threebusboundary}. The minimum damping ratio of the system decreases if the line admittances are increased and the system eventually becomes unstable. \label{fig:adm}}}
\end{figure}
Figure \ref{fig:adm} shows the minimum damping ratio $\zeta_{\min}$ as a function of the admittance of the line connecting inverter 1 to inverter 2 (blue), inverter 2 and 3 (red), and inverter 1 and 3 (orange). The black dashed line indicates the minimum damping ratio for the system in the original configuration described in Section \ref{sec.threebusboundary}.
In the original configuration, the admittance of the line from inverter 1 to inverter 2 is $0.0265$ S and the admittance of the lines connected to inverter 3 is $0.1327$ S. It can be seen that the system becomes increasingly underdamped and eventually unstable if any of the individual line admittances is increased beyond a threshold. This result is in line with the analytical insights obtained from Theorem \ref{thm.main} and confirms that adding, upgrading, or shortening transmission lines can make the system unstable.
Moreover, the conditions of Theorem~\ref{thm.main} and Proposition~\ref{cond.stab.power} depend on the maximum weighted node degree of the transmission network graph $d_{\max} = \max_{k \in \mc N} \sum\nolimits_{j:(j,k) \in \mc E} \norm{Y_{jk}}$, i.e., the maximum of the sum of the admittances of lines connected to individual inverters. This dependence can also be observed in Figure \ref{fig:adm}. Notably, the minimum damping ratio $\zeta_{\min}$ is insensitive to the admittance of the line connecting inverter 1 to inverter 2 (i.e., $\norm{Y_{12}}$) for small values of $\norm{Y_{12}}$, and becomes sensitive to $\norm{Y_{12}}$ only once $\norm{Y_{12}}$ becomes large enough to affect $d_{\max}$. 
}

\subsection{Stability boundaries of a three-bus transmission system}\label{sec.threebusboundary}
\begin{figure}[b!!!]
\input{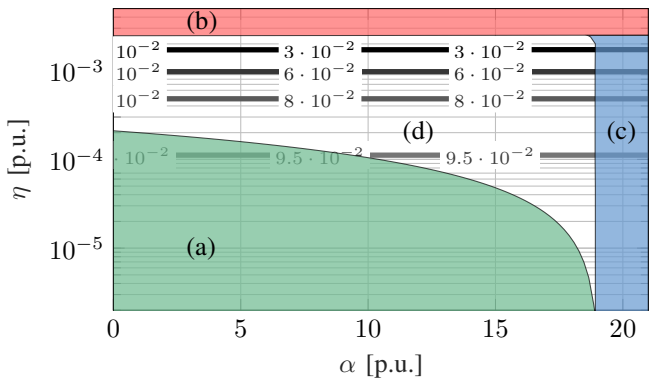}
\caption{Stability regions in parameter space. For control gains in region (a) Theorem \ref{thm.main} guarantees almost global asymptotic stability, in region (b) the system is unstable, in region (c) operational limits of power inverters are exceeded, and in region (d) the system is locally asymptotically stable. The lines in (d) indicate the minimum damping ratios of the linearized system.\label{fig.bound}}
\end{figure}

{We now consider the three-inverter transmission system shown in Figure \ref{fig:threebus} with the nominal line parameters and set-points given in Section~\ref{sec:three:adm}. We validate Theorem~\ref{thm.main} numerically by investigating the stability properties as a function of the control gains $\alpha$ and $\eta$.} The results are shown in Figure~\ref{fig.bound}. For control gains in region (a), Theorem~\ref{thm.main} guarantees almost global asymptotic stability, whereas for control gains in region (b) instability can be verified both via simulation or linearization. In region (c) the system remains stable in simulations of black starts and changes in load, but the magnitude $\norm{i_{o,k}}$ of the inverter output currents exhibits overshoots of more than $20 \%$. Due to tight limits on the maximum output current of power inverters, this is not desirable and, in practice, would require to oversize the inverters. Finally, for control gains in region (d), local asymptotic stability can be verified via linearization, and we observe that simulations of black starts and changes in load converge to $\mc T$. However, we cannot rule out the existence of unstable solutions in region (d), i.e., the union of (a) and (d) is an outer approximation of the range of parameters for which the system is almost globally asymptotically stable and satisfies the current limits of power inverters. Moreover, the lines in (d) indicate the minimum damping ratio {$\zeta_{\min}$} of the linearized system. In this example, minimum damping ratios below $5\cdot10^{-2}$ can result in significant oscillations and should be avoided. 

It should be noted that Condition~\ref{cond.stab} ensures exponential phase stability of the reduced-order system~\eqref{eq.voltage.control.kh}, i.e., using the same steps as in Prop.~\ref{thm.lyap.red.system} it can be shown that $\ddt \frac{1}{2}\norm{v}^2_{\mc S} \leq -\eta c \norm{v}^2_{\mc S}$. Thus, Theorem~\ref{thm.main} excludes region (c), and more generally, regions of the parameter space that result in poor damping. Thus, although the bound given by Theorem~\ref{thm.main} is conservative by an order of magnitude, the test case confirms that the controller gain must be limited to maintain stability, and it must be further limited to avoid oscillations and satisfy constraints of power inverters. In fact, within these operational constraints, the explicit bound given by Theorem~\ref{thm.main} is fairly accurate. We stress that Theorem~\ref{thm.main} is not restricted to operating points with zero power flow and gives almost global guarantees. Because of this, it is expected that the resulting bounds are conservative. 
\definecolor{mycolor1}{rgb}{0.00000,0.44700,0.74100}%
\definecolor{mycolor2}{rgb}{0.85000,0.32500,0.09800}%
\definecolor{mycolor3}{rgb}{0.92900,0.69400,0.12500}%

\subsection{Illustrative example: IEEE 9-bus system}\label{sec:example:ieee9}
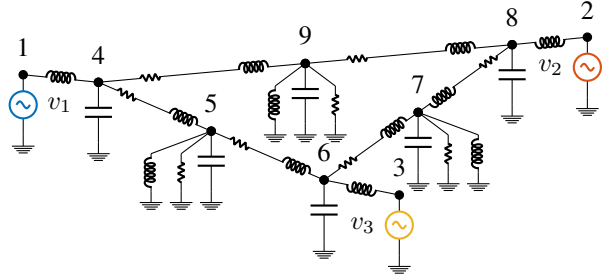
\begin{figure}[t!!!]
\begin{circuitikz}[american voltages]

\ctikzset{bipoles/resistor/height=0.15}
\ctikzset{bipoles/resistor/width=0.4}

\ctikzset{bipoles/generic/height=0.15}
\ctikzset{bipoles/generic/width=0.4}

\ctikzset{bipoles/length=.6cm}

\ctikzset{bipoles/thickness=2}

\node[label=1] at (0,0.1)  (I1)  {};
\node[label=2] (I2) at (7.5,0.6) {};
\node[label=3] (I3) at (5,-1.5) {};

\draw ($ (I1) + (0,-0.8) $) node [ground] (g1) {};
\draw (I1) to[sV, l=$v_1$, color=mycolor1] (g1);
\draw ($ (I2) + (0,-0.8) $) node [ground] (g2) {};
\draw (g2) to[sV, l=$v_2$, color=mycolor2] (I2);
\draw ($ (I3) + (0,-0.8) $) node [ground] (g3) {};
\draw (g3) to[sV, l=$v_3$, color=mycolor3] (I3);

\node[label=4] (4) at ($ (I1) + (1,-0.1) $) {};
\node[label=8] (8) at ($ (I2) + (-1,-0.1) $) {};
\node[label=6] (6) at ($ (I3) + (-1,0.2) $) {};

\draw ($ (4) + (0,-0.8) $) node [ground] (g4) {};
\draw (4) to[C] (g4);
\draw ($ (8) + (0,-0.8) $) node [ground] (g8) {};
\draw (8) to[C] (g8);
\draw ($ (6) + (0,-0.8) $) node [ground] (g6) {};
\draw (6) to[C] (g6);

\node[label=5] (5) at ($(4)!0.5!(6)$) {};
\node[label=9] (9) at ($(4)!0.5!(8)$) {};
\node[label=7] (7) at ($(6)!0.5!(8)$) {};

\draw ($ (5) + (0,-0.8) $) node [ground] (g5) {};
\draw (5) to[C] (g5);
\draw ($ (9) + (0,-0.8) $) node [ground] (g9) {};
\draw (9) to[C] (g9);
\draw ($ (7) + (0,-0.8) $) node [ground] (g7) {};
\draw (7) to[C] (g7);

\draw ($ (9) + (0.4,-0.8) $) node [ground] (g991) {};
\draw (g991) to[R] ($ (9) + (0.4,-0.3) $);
\draw (9) to[short] ($ (9) + (0.4,-0.3) $);
\draw ($ (9) + (-0.4,-0.8) $) node [ground] (g992) {};
\draw (g992) to[L] ($ (9) + (-0.4,-0.3) $);
\draw (9) to[short] ($ (9) + (0-0.4,-0.3) $);

\draw ($ (5) + (-0.4,-0.8) $) node [ground] (g551) {};
\draw (g551) to[R] ($ (5) + (-0.4,-0.3) $);
\draw (5) to[short] ($ (5) + (-0.4,-0.3) $);
\draw ($ (5) + (-0.8,-0.8) $) node [ground] (g552) {};
\draw (g552) to[L] ($ (5) + (-0.8,-0.3) $);
\draw (5) to[short] ($ (5) + (-0.8,-0.3) $);

\draw ($ (7) + (0.4,-0.8) $) node [ground] (g771) {};
\draw (g771) to[R] ($ (7) + (0.4,-0.3) $);
\draw (7) to[short] ($ (7) + (0.4,-0.3) $);
\draw ($ (7) + (0.8,-0.8) $) node [ground] (g772) {};
\draw (g772) to[L] ($ (7) + (0.8,-0.3) $);
\draw (7) to[short] ($ (7) + (0.8,-0.3) $);

\draw 		(I1) 
to[L,*-] (4);
\draw 		(I2) 
to[L,*-] (8);
\draw 		(I3) 
to[L,*-] (6);

\node (49) at ($(4)!0.5!(9)$) [label={[xshift=0cm, yshift=0.3cm]}]{};
\draw 		(4) 
to[R,*-] (49)
to [L,-*] (9);

\node (98) at ($(8)!0.5!(9)$) [label={[xshift=0cm, yshift=0.3cm]}]{};
\draw 		(9) 
to[R,*-] (98)
to [L,-*] (8);

\node (45) at ($(4)!0.5!(5)$) [label={[xshift=0cm, yshift=0.3cm]}]{};
\draw 		(4) 
to[R,*-] (45)
to [L,-*] (5);

\node (56) at ($(6)!0.5!(5)$) [label={[xshift=0cm, yshift=0.3cm]}]{};
\draw 		(5) 
to[R,*-] (56)
to [L,-*] (6);

\node (67) at ($(6)!0.5!(7)$) [label={[xshift=0cm, yshift=0.3cm]}]{};
\draw 		(6) 
to[R,*-] (67)
to [L,-*] (7);

\node (78) at ($(8)!0.5!(7)$) [label={[xshift=0cm, yshift=0.3cm]}]{};
\draw 		(8) 
to[R,*-] (78)
to [L,-*] (7);
\end{circuitikz}
\caption{IEEE 9-bus system. All generators have been replaced with inverters of the same rating that implement the control law~\eqref{eq.control.law}. The load dynamics are given by passive RLC circuits.\label{fig:ieee9}}
\end{figure}

In this section, we use the IEEE 9-bus system and replace generators by power inverters with the same rating. We use a structure preserving model with passive loads shown in Figure \ref{fig:ieee9} and do not modify the test-case parameters, i.e., Assumption~\ref{ass.constant.ratio} does not hold. 

\begin{figure}[t!!!]
\input{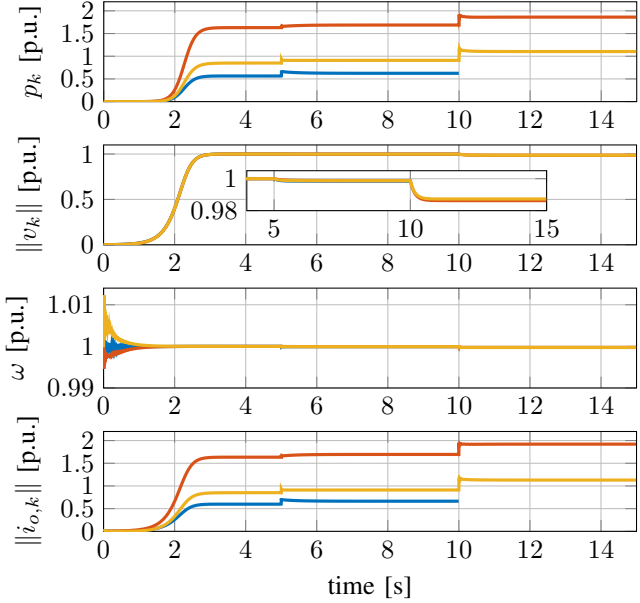}
\caption{Simulation of the IEEE 9-bus system. We show a grid black start at $t=0$s, a $20\%$ load increase at bus 5 at $t=5$s, and the loss of inverter $1$ at $t=10$s. The system is stable for a sufficiently small gain ($\eta = 10^{-3}$ p.u.).\label{fig.IEEE9}}
\end{figure}

\begin{figure}[h!!!]
\input{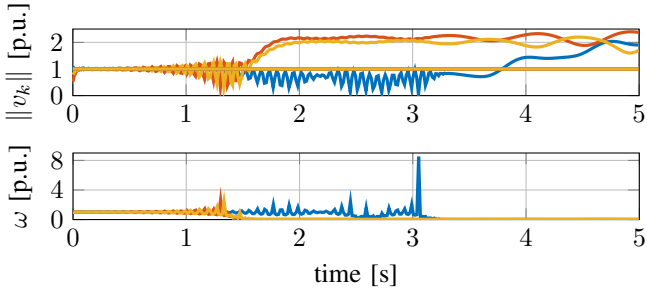}
\caption{Simulation of the IEEE 9-bus system for a large gain ($\eta = 10^{-2}$ p.u.). Instability is caused by high-gain control interfering with the line dynamics.\label{fig.IEEE9unstable}}
\end{figure}

We can therefore no longer apply Theorem~\ref{thm.main} directly. However, we illustrate how the dVOC~\eqref{eq.control.law} indeed synchronizes the grid to the desired power flows under nominal conditions and behaves well during contingencies. Furthermore, we illustrate {numerically} that, even if the assumptions of Theorem~\ref{thm.main} do not hold, the controller gain $\eta$ must remain small in order to ensure stability of the power network.
We use $\eta = 10^{-3}$ p.u., $\alpha = 10$ p.u., and the $\ell / r$ ratio of the lines is approximated by $\omega_0 \rho = 10$. We simulate the following events:  
\begin{itemize}
\item $t = 0\,s$ black start: $\|v_k(0)\|\approx 10^{-4}$ p.u.
\item $t = 5\,s$ $20\,\%$ active power increase at load $5$
\item $t = 10\,s$ loss of inverter $1$.
\end{itemize}
The results are illustrated in Figure~\ref{fig.IEEE9}. The controllers are capable to black-start the grid and converge to a synchronous solution with the desired power injections. When the load is increased ($t=5$\,s) we observe a droop-like behavior, the inverters maintain synchrony and share the power needed to supply the loads. Finally at $t=10$\,s we simulate a large contingency (the loss of Inverter 1). Inverters 2 and 3 do not loose synchrony and step up their power injection to supply the loads. 
Note that the current transients are particularly well behaved and do not present undesirable overshoots that are typical of other control strategies (e.g. machine emulation). In Figure~\ref{fig.IEEE9unstable} we show that, if we increase the gain $\eta$ to $10^{-2}$ the system is unstable. The fact that high gain control in conjunction with the line dynamics {is unstable is in agreement with the predictions made by Theorem~\ref{thm.main}}.

\section{Conclusion and outlook}\label{sec.conclusion}
In this paper, we considered the effect of the transmission-line dynamics on the stability of dVOC for grid-forming power inverters. A detailed stability analysis {for transmission lines with constant inductance to resistance ratio} was provided which shows that the transmission line dynamics have a destabilizing effect on the {multi-inverter} system, and that the gains of the inverter control need to be chosen appropriately. These instabilities cannot be detected using the standard quasi-steady-state approximation that is commonly used in power system stability analysis. Using tools from singular perturbation theory, we obtained explicit bounds on the controller gains and set-points that guarantee (almost) global asymptotic stability of the inverter based AC power system with respect to a synchronous steady state with the prescribed power-flows. {Broadly speaking, our conditions require that that the network is not to heavily loaded and that there is a sufficient time-scale separation between the inverter dynamics and line dynamics. Although the theoretical bounds are only sufficient and only apply for transmission lines with constant inductance to resistance ratio, we used a realistic test-cases to illustrate that the main salient features uncovered by our theoretical analysis translate to realistic scenarios, i.e., the power system becomes unstable when our sufficient stability conditions are violated.}
Similar instability phenomena induced by the dynamics of the transmission lines were recently observed in~\cite{vorobev2017high} for standard droop control. In view of these recent results, we believe that there is a need for more detailed studies to understand the fundamental limitations for the control of power inverters that arise from the dynamics of {transmission lines with heterogeneous inductance to resistance ratios, transformers, and other network dynamics that are typically not considered in power system stability analysis.}

\appendix 
{
\emph{Proof of Proposition \ref{prop.sync}:}
 We can express the phase error as
\begin{align*}
 e_{\theta,k} 
 \!=\!\!\sum_{(j,k)\in\mc E} \norm{Y_{jk}}\! (v_j\!-v_k) + \!\!\sum_{(j,k)\in\mc E}  \!\norm{Y_{jk}}{\left(I_2\!-\!\tfrac{v_j^\star}{v_k^\star}R(\theta_{jk}^\star)\right)} v_k.
\end{align*}
Using Assumption 1, it can be verified that $\sum_{(j,k)\in\mc E} \norm{Y_{jk}} (v_j\!-v_k) = -\mc R(\kappa) i^s_o(v)$. Next, we use
$\sin(\kappa) =  \tfrac{\ell_{jk}\,\omega_0}{\sqrt{r_{jk}^2+\omega_0^2\ell_{jk}^2}}$ and $\cos(\kappa)  = \tfrac{r_{jk}}{\sqrt{r_{jk}^2+\omega_0^2\ell_{jk}^2}}$
to write $p^\star_{jk}$ and $q^\star_{jk}$ (cf. Condition \ref{cond.consistent}) as
\begin{align}
\begin{split}\label{eq.powers.ref}
p^\star_{jk} ={v_k^{\star2}} \norm{Y_{jk}} \left( \cos(\kappa)  -\tfrac{v^\star_j}{v^\star_k} \cos(\theta_{jk}^\star -\kappa)\right),\\
 q^\star_{jk} = {v_k^{\star2} } \norm{Y_{jk}} \left(  \sin(\kappa) +\tfrac{v^\star_j}{v^\star_k}  \sin(\theta_{jk}^\star -\kappa)\right).
\end{split}
\end{align}
Using $p^\star_k = \sum_{(j,k)\in\mc E} p^\star_{jk}$, $q^\star_k = \sum_{(j,k)\in\mc E} q^\star_{jk}$, and~\eqref{eq.powers.ref} we obtain
\begin{align}\label{eq.powang}
\!\!\begin{bmatrix}
p_k^\star\! & q_k^\star\\
-q_k^\star\! & p_k^\star
\end{bmatrix}\! = 
v_k^{\star2} R(\kappa)^{\mathsf{T}} \!\! \sum_{(j,k)\in\mc E} \! \norm{Y_{jk}}\left(I_2-\tfrac{v^\star_j}{v^\star_k}R(\theta_{jk}^\star)\right)\!.\!
\end{align}
and the proposition immediately follows.\hfill\IEEEQED

\emph{Proof of Proposition \ref{cond.stab.power}:} 
We begin by observing that 
\begin{align*}
\left|1-\tfrac{v^\star_j}{v^\star_k}\cos(\theta^\star_{jk})\right| = \left|\sin(\kappa)^2 + \cos(\kappa)^2-\tfrac{v^\star_j}{v^\star_k}\cos(\theta^\star_{jk})\right|.
\end{align*}
Using the identity $\cos(\theta)=\cos(a+b) = \cos(a)\cos(b)-\sin(a)\sin(b)$ with $a=\theta-\kappa, b=\kappa$, applying the triangle inequality, and noting that $0 \leq \kappa \leq \tfrac{\pi}{2}$ results in
\begin{align*}
\left|1-\tfrac{v^\star_j}{v^\star_k}\cos(\theta^\star_{jk})\right| \leq &\cos(\kappa) \left|\cos(\kappa)-\tfrac{v^\star_j}{v^\star_k}\cos(\theta-\kappa)\right|\\
&+\sin(\kappa) \left| \sin(\kappa) + \tfrac{v^\star_j}{v^\star_k}\sin(\theta-\kappa)\right|.
\end{align*}
Using \eqref{eq.powers.ref} it follows that $\norm{Y_{jk}}\left|1-\tfrac{v^\star_j}{v^\star_k}\cos(\theta^\star_{jk})\right| \leq \tfrac{\cos(\kappa)}{v^{\star2}_k} \left|p^\star_{jk}\right| +\tfrac{\sin(\kappa)}{v^{\star2}_k} \left| q^\star_{jk} \right|$. Substituting into Condition~\ref{cond.stab} and letting $\bar{\theta}^\star=\tfrac{\pi}{2}$ completes the first part of the proof. Next, we note that $\norm{K_k} = \frac{1}{v^{\star2}_k} \sqrt{p^{\star2}_k + q^{\star2}_k}$, i.e., $\norm{K_k} = s^\star_k {v^\star_k}^{-2}$, and it follows that $\norm{\mc K} = \max_{k \in \mc N} s^\star_k {v^\star_k}^{-2}$. Moreover, given matrices $F_1 \in \mathbb{R}^{n_1 \times m_1}$ and $F_2 \in \mathbb{R}^{n_2 \times m_2}$  it holds that $\norm{F_1 \otimes F_2} = \norm{F_1} \norm{F_2}$ (see \cite[p. 413]{lancaster1972norms}). Using this fact, it holds that $\norm{\mc L} = \norm{L} \leq 2 \max_k \sum_{j:(j,k) \in \mc E} \norm{Y_{jk}}$. Finally, because $\mc R(\kappa)$ is a rotation matrix it holds that $\norm{\mc Y} = \norm{\mc R(\kappa)^\mathsf{T} \mc L} = \norm{L}$ and the proposition follows. 
\hfill\IEEEQED
}
 
\emph{Proof of Theorem~\ref{thm.AGAS}:}
 We first establish stability of $\mc C$. It follows from $\chi_3 \in \mathscr{K}$ that $\ddt \mc V \leq 0$. This implies that $\chi_1(\norm{\varphi_f(t,x_0)}_{\mc C}) \leq \mc V(\varphi_f(t,x_0)) \leq \mc V(x_0) \leq \chi_2(\norm{x_0}_{\mathcal{C}})$. Next, we use $\chi^{-1}_1 \in \mathscr{K}_{\infty}$ denote the inverse function of $\chi_1 \in \mathscr{K}_{\infty}$ and obtain $\norm{\varphi_f(t,x_0)}_{\mc C} \leq \chi^{-1}_1(\chi_2(\norm{x_0}_{\mc C}))$ for all $t \in \R_{\geq 0}$. It follows from standard arguments (see \cite[Sec. 25]{H67}) that $\mc C$ is Lyapunov stable according to Definition~\ref{def:ags}. Next, we establish almost global attractivity of $\mc C$. Because $\mc U$ is invariant with respect to $\ddt x = f(x)$, points in $\mc U$ are only reachable from initial conditions $x_0 \in \mc Z_{\mc U} \supseteq \mc U$. Thus, for $x_0 \notin \{\mc Z_{\mc U} \cup \mc C\}$ and for all $t \in \mathbb{R}_{>0}$ such that $\varphi_f(t,x_0) \notin \mc C$ it holds that $\ddt \mc V(\varphi_f(t,x_0)) < 0$. Using standard arguments (see \cite[Sec. 25]{H67}) and noting that $\mc V(x)=0$ for all $x\in\mc C$ it follows that $\lim_{t \to \infty} \mc V(\varphi_f(t,x_0)) \to 0$. Moreover, using $\norm{x}_{\mc C} \leq \chi^{-1}_1(\mc V(x))$ we conclude that 
\begin{align}
 \forall x_0 \notin \mc Z_{\mc U}: \lim_{t \to \infty} \norm{\varphi_f(t,x_0)}_{\mc C} = 0.
\end{align}
Because $\mc Z_{\mc U}$ has zero Lebesgue measure this is precisely the definition of almost global attractivity given in Definition~\ref{def:ags}. Because $\mathcal{C}$ is almost globally attractive and stable, it is almost globally asymptotically stable.\hfill\IEEEQED

\begin{lemma}\label{lem.pos}
Let $v_1^\star,...,v_N^\star\in\R_{>0}$. The polynomials 
\begin{equation}\label{eq.positive}
p_m(x) = \sum_{k=1}^N \bigg(\frac{\sum_{n=1}^Nv_n^{\star2}}{{v^{\star{m-1}}_k}}x_k^{m+1} -  \sum_{j=1}^N\frac{v^\star_k v^\star_j}{v^{\star m-1}_k}x_k^{m}x_j\bigg) 
\end{equation}
are nonnegative on the nonnegative orthant for all $m\in\mathbb N$.
\end{lemma}
\begin{IEEEproof}
We prove the Lemma by induction. Letting 
$
s \coloneqq \begin{bmatrix}
v^\star_1&
\cdots&
v^\star_N
\end{bmatrix}^\mathsf{T}\!\!
$, $p_1(x)$ can be written as $p_1(x)=x^\mathsf{T} \left(\sum\nolimits_{n=1}^Nv_n^{\star2} I_N - ss^\mathsf{T}\right)x\ge 0$. 
Now we define $\bar x \coloneqq \frac{1}{\sum_{n=1}^Nv_n^{\star2}}\sum_{j=1}^N v^\star_j{x_j}$. Note that $\bar x \ge0$ whenever $x\in\R^N_{\ge0}$.
We now assume that $m\ge 2$ and obtain
\begin{align*}
p_m(x) 
&=   \sum_{k=1}^N \left(\frac{\sum_{n=1}^Nv_n^{\star2}}{{v^{\star m-1}_k}}x_k^m \left( x_k -{v_k^\star}\, \bar x \right)\right)\\
&=  \bar x \sum_{k=1}^N \left(\frac{\sum_{n=1}^Nv_n^{\star2}}{{v^{\star m-2}_k}}x_k^{m-1} \left( x_k -{v_k^\star}\, \bar x \right)\right)  \\
&+ \sum_{k=1}^N \frac{\sum_{n=1}^Nv_n^{\star2}}{{v^{\star m-1}_k}}x_k^{m-1} \left( x_k -{v_k^\star}\, \bar x \right)^2\\
& = \bar x p_{m-1}(x)+ \sum_{k=1}^N\frac{\sum_{n=1}^Nv_n^{\star2}}{{v^{\star m-1}_k}}x_k^{m-1}  \left( x_k -{v_k^\star}\, \bar x \right)^2.
\end{align*}
If $p_{m-1}(x)$ is nonnegative on the nonnegative orthant it follows that $p_m(x)\ge0$ for all $x\in\R^N_{\ge0}$ and the proof is complete.
\end{IEEEproof}

\emph{Proof of Lemma \ref{lem.positivsatz.1}:}
Since
\begin{align*}
v^\mathsf{T} P_S \Phi(v) v  =  v^\mathsf{T} P_S v -  v^\mathsf{T}  \diag\left(\left\{\tfrac{\norm{v_k}^2}{v_k^{\star2}} I_2\right\}_{k=1}^N\right) P_S v, 
\end{align*}
the inequality~\eqref{eq.claim1} is equivalent to  
\begin{align*}
  v^\mathsf{T}  \diag\left(\left\{\tfrac{\norm{v_k}^2}{v_k^{\star2}} I_2\right\}_{k=1}^N\right) P_S v \geq 0, \quad \forall v\in\R^{n}. 
\end{align*}
Next, we use $\theta_{kj}$ to denote the relative angle between $v_j$ and $v_k$ such that $\tfrac{v_j}{\norm{v_j}} = R(\theta_{kj})\tfrac{v_k}{\norm{v_k}}$ holds. Given the particular form of $P_S$ we can write  
\begin{align*}
& \sum_{n=1}^Nv_n^{\star2}~ v^\mathsf{T}  \diag\left(\left\{\frac{1}{v_k^{\star2}}\norm{v_k}^2 I_2\right\}_{k=1}^N\right) P_S v\\
& = \sum_{k=1}^N \bigg(\frac{\sum_{n=1}^Nv_n^{\star2}}{{v^{\star2}_k}}\|v_k\|^4 -  \sum_{j=1}^N\frac{ v^\star_j}{v^\star_k}\|v_k\|^3\norm{v_j}\cos(\theta_{kj} - \theta^\star_{kj})\bigg) \\
&\ge \sum_{k=1}^N \bigg(\frac{\sum_{n=1}^Nv_n^{\star2}}{{v^{\star2}_k}}\|v_k\|^4 -  \sum_{j=1}^N\frac{ v^\star_j}{v^\star_k}\|v_k\|^3\norm{v_j}\bigg)\\
& = p_3\left(\,[\norm{v_1},...,\norm{v_N}]^\mathsf{T}\,\right)\ge0,
\end{align*}
where $p_3(\cdot)$ is defined in~\eqref{eq.positive}, and the last inequality follows from Lemma~\ref{lem.pos}. It follows that the inequality \eqref{eq.claim1} holds.
\hfill\IEEEQED

\emph{Proof of Lemma~\ref{lem.redlfdecr}:~}
 We start by noting that, since $(\mc K-\mc L)v=0$ for all $v\in\mc S$ and $P_S$ is the projector onto $\mc S ^\perp$, it holds that $(\mc K-\mc L) = (\mc K-\mc L)P_S$.
Thus,~\eqref{eq.ass.reformulation} is equivalent to
\begin{align}
\! v^{\!\mathsf{T}} \!(P_S \mc K P_S \!+\!\alpha P_S) v \!\leq\! v^{\!\mathsf{T}}\! P_S \mc L P_S v - c\, \|v\|_{\mc S}^2. \label{eq.ineq.condition}
\end{align}
Noting that 
{$\frac{1}{2}(K_k + K^\mathsf{T}_k) = \sum_{j=1}^N \norm{Y_{jk}}(1- \tfrac{v^\star_j}{v^\star_k}\cos(\theta^\star_{jk}))$ for all $k\in\mc N$ (cf. Proposition \ref{prop.sync}), 
the left-hand side of~\eqref{eq.ineq.condition} can be bounded as follows
\begin{align}
&v^\mathsf{T} [P_S\mc K P_S +\alpha P_S] v \notag \\ &\; \le  \|v\|_{\mc S}^2\left[\max_k \sum\nolimits_{j=1}^N \norm{Y_{jk}}\left|1-\tfrac{v^\star_j}{v^\star_k}\cos(\theta^\star_{jk})\right|+\alpha\right]\!. \label{eq.bound3}
\end{align}
After some lengthy algebraic manipulations (see \cite[Prop. 7]{colombino2017global2}) the quadratic form $v^\mathsf{T} P_S \mc L P_S v$ can be bounded as follows
  \begin{align*}
 v^\mathsf{T} P_S \mc L P_S v       \ge \lambda_2( L)\|v\|_{\mc S}^2\frac{1}{N{\sum_{n=1}^Nv^{\star2}_n}}  \sum_{k=1}^N\sum_{j=1}^N v^\star_j{v^\star_k} \cos(\theta^\star_{jk})
\end{align*}
Next, let $\lceil x \rceil\coloneqq \min_{y\in\mathbb N, y\ge x} y$ denote the ceiling operator that rounds up a real number $x \in \mathbb{R}^n$ to the nearest integer. The fact that $\sum_{k=1}^N\sum_{j=1}^Nv^\star_j{v^\star_k}\cos(\theta^\star_{jk})\ge N^2 v^{\star2}_{\min}\frac{1}{2}\left(1{+}\cos(\bar{\theta}^\star)\right)$ can be proven by noting that a minimizer for the sum of cosines is given by $\theta^\star_{k1} = 0$ for all $k \in \{2,\ldots,\lceil\frac{N}{2}\rceil\}$ and $\theta^\star_{k1}=
\bar{\theta}^\star$ for all $k \in \{\lceil\frac{N}{2}\rceil+1,\ldots,N\}$. It immediately follows that
\begin{align}
   v^\mathsf{T} P_S \mc L P_S v  \ge \frac{1}{2}\frac{v^{\star2}_{\min}}{v^{\star2}_{\max}}\left(1 {+} \cos(\bar{\theta}^\star)\right) \lambda_2( L)\|v\|_{\mc S}^2. \label{eq.bound2}
\end{align}
By substituting~\eqref{eq.bound3} and~\eqref{eq.bound2} into~\eqref{eq.ineq.condition}, we can conclude that Condition~\ref{cond.stab} is a sufficient condition for~\eqref{eq.ass.reformulation}.\hfill\IEEEQED
%

\emph{Proof of Proposition~\ref{prop.beta1}:}
The function $f_v(v,i)$ in~\eqref{eq.voltage.cl.rot} is separable in its two arguments and linear in $i$. Hence,
\begin{align}\label{eq.beta1.def}
\!\!f_v(v,y+i^s(v)) - f_v(v,i^s(v)) = f_v(0,y) = \eta {\mc R(\kappa)}  \mc B y.\!
\end{align}
Using $y_o = \mc B y$, it directly follows that 
\begin{align*}
&\frac{\partial V}{\partial v}[f_v(v,y+i^s(v))-f_v(v,i^s(v))]\\
&\hphantom{\frac{\partial V}{\partial v}}\le (\eta v^\mathsf{T} P_S + 2 \alpha \alpha_1 \eta^2 v^\mathsf{T} \Phi(v)){\mc R(\kappa)}y_o\\
&\hphantom{\frac{\partial V}{\partial v}}\le (\eta \|v\|_{\mc S}  + 2 \alpha \alpha_1 \eta^2 \norm{\phi(v) v})\norm{y_{o}} \le  \beta_{yv} \psi(v) \norm{y_o},
\end{align*}
with 
$
\beta_{yv}  \coloneqq \max\left(\norm{\mc K -\mc L}^{-1}, {2 \alpha_1 \eta} \right).
$
Next, using~\eqref{eq.ass.reformulation} it can be verified that
$
c \,\|v\|_{\mc S}^2  +v^\mathsf{T} \alpha I_{2N} v \leq  - v^\mathsf{T} P_S ( \mc K - \mc L)P_Sv
$. This results in 
$ 
c \,\|v\|_{\mc S}^2   \leq  \|v\|_{\mc S}^2\| \mc K - \mc L\|$ and $c\le \norm{\mc K -\mc L}$. From~\eqref{thm1.strictbound} and the fact that $c\le \norm{\mc K -\mc L}$ we conclude that
$
2 { \alpha_1 \eta } \le \frac{2 c}{5 \norm{\mc K -\mc L}^2} \le \beta_1 = \norm{\mc K -\mc L}^{-1}
$
and the proof is complete. \hfill\IEEEQED

\emph{Proof of Proposition~\ref{prop.beta2}:}
Using~\eqref{eq.beta1.def} and $y_o = \mc B y$, we obtain
\begin{align*}
&-\frac{\partial W}{\partial y}\frac{\partial i^s}{\partial v} f_v(v,y+i^s(v)) = \rho y^\mathsf{T}_o \mc Y (f_v(v,i^s(v)) + f_v(0,y))\\
&\hphantom{-\frac{\partial W}{\partial y}} + \rho y^\mathsf{T} L_T \mc B_n \mc B^\mathsf{T}_n L_T Z^{-1}_T \mc B^\mathsf{T}(f_v(v,i^s(v)) + f_v(0,y)).
\end{align*}
Moreover, it holds that $L_T Z^{-1}_T =\frac{\rho}{\omega_0^2 \rho^2 +1} (I_{2M}- \omega_0 \mc J_M)$, and it follows from the mixed-product property of the Kronecker product
that $\mc B^\mathsf{T}_n \mc B^\mathsf{T}=(B B_n)^\mathsf{T} \otimes I_2 = \mathbbl{0}_{2M}$, and $\mc B^\mathsf{T}_n \mc J_M \mc B^\mathsf{T} = (B B_n)^\mathsf{T} \otimes J =  \mathbbl{0}_{2M}$. Thus, it holds that $\mc B^\mathsf{T}_n L_T Z^{-1}_T \mc B^\mathsf{T} = \mathbbl{0}_{2M}$ and the proposition immediately follows from
$\|f_v(v,i^s(v))\| \leq \psi(v)$ and $\norm{f_v(0,y)} = \eta \norm{y_o}$.\hfill\IEEEQED
\bibliographystyle{IEEEtran}
\bibliography{IEEEabrv,bib_file}

\end{document}